\documentclass[11pt,a4paper]{article}

\usepackage{graphicx}
\usepackage{amsmath}
\usepackage{amsfonts}
\usepackage{amssymb}
\usepackage{enumerate}
\usepackage{lscape}
\usepackage{longtable}
\usepackage{rotating}
\usepackage{color}
\usepackage{url}
\usepackage{rotating}
\usepackage{algpseudocode}
\usepackage[ruled]{algorithm}
\usepackage{pgf,tikz}
\usepackage{subfig}
\usepackage[titletoc]{appendix}

\usepackage{mathrsfs}
\usepackage{pgfplots}
\usetikzlibrary{arrows}
\usepackage{color}

\usepackage{bbm}

\usepackage{booktabs}
\usepackage{multirow}
\usepackage{mathtools}
\usepackage{makecell}
\usepackage[titletoc]{appendix}

\usepackage[round]{natbib}
\usepackage{eqparbox} 

\newcommand{\mathbox}[3][\mathop]{
	#1{\eqmakebox[#2]{$\displaystyle#3$}}%
}

\numberwithin{equation}{section} 

\usepackage[symbol]{footmisc} 

\DeclareMathOperator*{\argmin}{argmin}

\DeclareMathOperator*{\sub}{sub}
\DeclareMathOperator*{\maxlim}{maxlim}
\DeclareMathOperator*{\minlim}{minlim}

\DeclareMathOperator*{\emis}{emis}
\DeclareMathOperator*{\dep}{dep}
\DeclareMathOperator*{\DEM}{DEM}
\DeclareMathOperator*{\G}{G}
\DeclareMathOperator*{\V}{V}
\DeclareMathOperator*{\C}{C}

\DeclareMathOperator*{\N}{N}
\DeclareMathOperator*{\E}{E}

\DeclareMathOperator*{\SC}{SC}

\DeclareMathOperator*{\arr}{arr}

\DeclareMathOperator{\veh}{VEH}
\DeclareMathOperator{\rep}{rep}
\DeclareMathOperator{\cop}{cop}
\DeclareMathOperator{\TS}{TS}
\DeclareMathOperator{\CO}{CO}
\DeclareMathOperator{\capac}{cap}

\DeclareMathOperator{\TGE}{TGE}
\DeclareMathOperator{\TSC}{TSC}
\DeclareMathOperator{\TDT}{TDT}
\DeclareMathOperator{\MDD}{MDD}
\DeclareMathOperator{\LAC}{LAC}
\DeclareMathOperator{\MAT}{MAT}

\DeclareMathOperator{\ordy}{ord}

\DeclareMathOperator{\SP}{SP}

\DeclareMathOperator{\PP}{P}

\newcommand{\beqn}[1]{\begin{equation}\label{#1}}
\newcommand{\eeqn}{\end{equation}}

\usepackage{accents}
\newlength{\dhatheight}
\newcommand{\doublehat}[1]{%
    \settoheight{\dhatheight}{\ensuremath{\hat{#1}}}%
    \addtolength{\dhatheight}{-0.35ex}%
    \hat{\vphantom{\rule{1pt}{\dhatheight}}%
    \smash{\hat{#1}}}}

\definecolor{aau2}{rgb}{0.0, 0.5, 0.69}
\definecolor{aau3}{rgb}{0.0, 0.53, 0.74}
\definecolor{aau4}{rgb}{0.0, 0.48, 0.65}
\definecolor{aau5}{rgb}{0.0, 0.45, 0.73}
\definecolor{rsap}{RGB}{130, 36, 51}
\definecolor{gsap}{RGB}{112, 164, 137}

\definecolor{tud}{rgb}{0.43,0.73,0.11}
\definecolor{verde}{rgb}{0.33,0.53,0.11}
\definecolor{darkgreen}{rgb}{0,0.6,0}

\definecolor{ttffqq}{rgb}{0.0, 0.48, 0.65} 
\definecolor{ffqqqq}{rgb}{0.0, 0.5, 0.69} 

\usetikzlibrary{arrows}

\allowdisplaybreaks[1]

\tikzstyle{decision} = [diamond, draw, fill=blue!20,
text width=4.5em, text badly centered, node distance=3cm, inner sep=0pt]
\tikzstyle{block} = [rectangle, draw, fill=blue!20,
text centered, rounded corners, minimum height=4em]
\tikzstyle{line} = [draw, -latex']
\tikzstyle{cloud} = [draw, ellipse,fill=red!20, node distance=3cm,
minimum height=2em]
\tikzstyle{cloud2} = [draw, ellipse,fill=green!20, node distance=3cm,
minimum height=2em]
\begin{document}
	
	\title{An integrated assignment, routing, and speed model \\ for roadway mobility and transportation \\ with environmental, efficiency, and service goals}
	
	\author{
		T. Giovannelli\thanks{Department of Industrial and Systems Engineering, Lehigh University, Bethlehem, PA 18015, USA ({\tt tog220@lehigh.edu}).}
		\and
		L. N. Vicente\thanks{Department of Industrial and Systems Engineering, Lehigh University, Bethlehem, PA 18015, USA ({\tt lnv@lehigh.edu}). Support for this author was partially provided by the Centre for Mathematics of the University of Coimbra under grant FCT/MCTES UIDB/MAT/00324/2020.}
	}
	
	\maketitle
	
	\begin{abstract}
	
	    Managing all the mobility and transportation services with autonomous vehicles for users of a smart city requires determining the assignment of the vehicles to the users and their routing in conjunction with their speed. Such decisions must ensure low emission, efficiency, and high service quality by also considering the impact on traffic congestion caused by other vehicles in the transportation network. 
	    
	    In this paper, we first propose an abstract trilevel multi-objective formulation architecture to model all vehicle routing problems with assignment, routing, and speed decision variables and conflicting objective functions. Such an architecture guides the development of subproblems, relaxations, and solution methods. We also propose a way of integrating the various urban transportation services by introducing a constraint on the speed variables that takes into account the traffic volume caused across the different services.
	    Based on the formulation architecture, we introduce a (bilevel) problem where assignment and routing are at the upper level and speed is at the lower level. To address the challenge of dealing with routing problems on urban road networks, we develop an algorithm that alternates between the assignment-routing problem on an auxiliary complete graph and the speed optimization problem on the original non-complete graph. The computational experiments show the effectiveness of the proposed approach in determining approximate Pareto fronts among the conflicting objectives.
	\end{abstract}

	\section{Introduction} \label{sec:intro}

    In this paper, we propose an innovative integrated transportation model for the management of possibly all vehicles traveling on the streets and roads of a city, which are assumed to have different levels of autonomy (with or without a driver). Our main goal is to provide an optimization model that can be used to effectively manage mobility and transportation within a city by adopting a green logistics perspective and pursuing efficiency and service quality. 
    The integrated model introduced in this work will qualify the cities implementing the proposed transportation network as smart. 
    
    The achievement of our goal will be enabled by the recent (and future) advances in information and location-sensitive technologies, which facilitate acquiring the necessary high-quality (granular) data for supporting the decision-making process for vehicles traveling on the streets and roads of a smart city.
	Several definitions of smart city have been proposed in the literature (see, e.g.,~\cite{TNam_TAPardo_2011}, \cite{ZBronstein_2009}, and \cite{RWenge_XZhang_CDave_LChao_SHao_2014}) and a universal characterization cannot be provided since smart cities involve different types of features that vary based on the specific context considered~(\cite{AKarvonen_FCugurullo_FCaprotti_2018}).
	In general, the phrase {\it smart city} refers to cities provided with technological infrastructure based on advanced data processing that pursue several goals, such as a more efficient city governance, a better life quality for citizens, increasing economic success for businesses, and a more sustainable environment~(\cite{CYin_ZXiong_HChen_2015}).
	
	In our paper, we focus on smart cities with intelligent transportation systems~(see~\cite{ZXiong_2012} for a survey). In particular, all the smart cities provided with technologies for real-time transportation data acquisition fit within the scope of our work. However, we point out that the process leading to the acquisition of such technologies, which requires a multitude of social, political, and economic factors involving public-private partnerships and local authorities, is omitted from this paper. 
	The importance of the role played by mobility and transportation in smart cities is highlighted in many works (\cite{CYin_ZXiong_HChen_2015}). For instance, in~\cite{ZTang_KJayakar_XFeng_etal_2019}, four different groups of smart cities are identified based on a cluster analysis of cities around the world, and one of them gathers cities adopting smart transportation systems. The goal of such systems is to control traffic congestion by public transportation, car sharing, and self-driving cars. In~\cite{AArroub_BZahi_ESabir_etal_2016}, the urban congestion is seen as a challenge arising from the persistent need of citizens to use their private cars. A solution proposed to alleviate such a problem requires smart control of the traffic in the existing road infrastructure to ensure a sustainable transportation, which is exactly the goal of our paper.
	
	A crucial tool for managing transportation in smart cities is the concept of urban artificial intelligence (AI) (\cite{FCugurullo_2020}). Among the examples of urban AI, self-driving cars and city brains represent two important categories. In particular, the author points out that the number of cities where autonomous cars are allowed to drive is increasing~(see also~\cite{RAAcheampong_FCugurullo_2019}), and these also include cars with the highest level of autonomy, when no human input or supervision are required. A city brain is a digital platform applied to the management of a city, including urban transportation, where the goal is to  control traffic lights and flows of vehicles by using advanced data collected throughout the city. In this paper, we assume that our integrated model is used by a city brain for the urban transportation management of autonomous vehicles. However, the development of the features of such a platform is left for future work.
	
	The presence of autonomous vehicles in a transportation network allows for the determination of the assignment of vehicles to users and their routing in conjunction with their speed, thus increasing the complexity of the decision-making process, which is naturally formulated as a hierarchy of three levels (assignment, routing, and speed). Transportation and mobility solutions should be determined with low environmental impact but, at the same time, with high efficiency for service providers and high service quality for users, leading to conflicting goals which need to be optimized simultaneously. Moreover, decisions unilaterally made in one component or part of the network can have a strong impact on the overall system. Consequently, the optimal solutions of the single components considered separately are different from the solutions obtained when such interaction is taken into account, thus requiring a proper integration of the network components. Finally, the current modeling techniques for assignment-routing problems are based on assumptions that are not satisfied on urban road networks and, therefore, require innovative solution methodologies. We aim at developing a future-oriented \textit{integrated assignment, routing, and speed system for roadway mobility and transportation with environmental, efficiency, and service goals}, and this paper represents an instrumental building block towards this main goal and where we make the following three main contributions.
	
	\subsection*{A trilevel multi-objective formulation architecture for vehicle routing problems with speed optimization}
	We propose an innovative formulation architecture to decompose every \textit{Vehicle Routing Problem}~(VRP) with speed optimization into three levels, each addressing one of the following vehicle-related decisions: assignment to users, routing to accommodate all the requests, and speed along each segment of the route\footnote{Note that although controlling the speed in roadway transportation is considered an impractical task (see, e.g., \cite{TVidal_2020}), considering smart and autonomous vehicles allows one to take into account speed decisions.}. In particular, our paper extends the bilevel formulation developed by~\cite{YMarinakis_AMigdalas_PPardalos_2007} to propose a trilevel multi-objective formulation for a~VRP with speed optimization (referred to as a~VRP/speed problem). Such an architecture allows a comprehensive understanding of the overall problem complexity, guides the development of subproblems, relaxations, and solution methods, and provides new insights into~VRP/speed problems. In this paper, we use the trilevel multi-objective architecture to formulate a (bilevel) VRP/speed problem (where upper level is assignment-routing and lower level is speed), develop a corresponding optimization method, and identify and report the trade-offs among the three considered goals.
	
	\subsection*{Integration among all the different transportation problem components in a smart city}
    The transportation services arising in an urban transportation network can be modeled through a proper~VRP variant. All the frameworks, models, and methods proposed in the literature address the different transportation services arising in a city independently, without any attempt to consider such services as parts of the same system. Several studies in the literature considered modeling a general framework accounting for as many transportation services as possible. For example, \cite{TVidal_TCrainic_MGendreau_CPrins_2014} proposed a unified model capable of separately describing different transportation {\it problem components} (later referred to as {\it components}). However, the goal of~\cite{TVidal_TCrainic_MGendreau_CPrins_2014} was to propose a general-purpose optimization algorithm that can quickly provide efficient solutions to different problems, each related to a transportation service. In contrast, our paper aims to develop a new general model, where different transportation components (such as personal trips, freight transportation, ride-sharing, car-pooling, dial-a-ride, and vehicle sharing) can be integrated into a comprehensive optimization framework. In this regard, we propose a~VRP/speed formulation for a specific problem component, which is sufficient to fully address the integration among all the components. In particular, in such a formulation, an innovative constraint on the speed variables is used to model the impact of the traffic congestion caused by routing decisions (made in other components) on the component under consideration.
    
	\subsection*{Use of non-complete graphs for vehicle routing problems with speed optimization}
    Using non-complete graphs is essential to address a VRP/speed problem on urban road networks. However, to find the value of the speed variables for each edge in the network, we cannot directly resort to commonly used approaches for non-complete graphs (which build a complete graph by combining the edges in the shortest path between two customer nodes into a single edge). Therefore, we propose an approach that computes an approximate solution to the assignment-routing (upper level) problem on the road network by first solving such a problem on an auxiliary complete graph. 
    This mechanism allows one to use existing~VRP methods for complete graphs, without the need to formulate the assignment-routing problem on the non-complete graph, where some of the classical~VRP constraints are most likely infeasible. Then, considering the approximate assignment-routing solution as a parameter, we develop a formulation for the speed optimization (lower level) problem on a non-complete graph.  
    It is important to remark that VRP with green-oriented objectives and speed optimization is well-known in the literature~(\cite{TBekta_2011}). However, considering such a problem on a road network represents, to the best of our knowledge, a novel contribution, thus requiring a new solution methodology.
	
\subsection{Organization of this paper}

This paper is organized as follows. Section~\ref{sec:background_operations_research_models} provides a general overview of the main optimization models adopted to solve transportation problems on complete graphs and road networks. The trilevel multi-objective formulation architecture is described in Section~\ref{sec:formulation}. Section~\ref{sec:mathematical_programming} presents an innovative constraint to model the integration among all the different transportation components. In Section~\ref{sec:VRP_complete_graph}, the trilevel multi-objective architecture is used to formulate a (bilevel)~VRP/speed problem on a non-complete graph and an optimization method is developed to solve such a problem. The  computational experiments are described in Section~\ref{sec:computational_experiments}. Finally, in Section~\ref{sec:conclusion} we draw some concluding remarks and we outline several ideas for future work.
	
	\section{Literature review} \label{sec:background_operations_research_models}

	\subsection{Optimization models for transportation on complete graphs} \label{subsec:background_operations_research_models}

	Transportation services across a city have been widely studied within the field of Operations Research~(see, e.g.,~\cite{GKim_2015} and~\cite{DCattaruzza_2017}). Optimization problems related to such services are variants of the basic VRP, where a fleet of vehicles is used to deliver goods to a set of customer nodes. Two important decisions are considered: assigning groups of customers to each vehicle and defining the corresponding route.
	
	\textit{The Pickup and Delivery Problem} (PDP) is a VRP variant where people or objects need to be transported from an origin to a destination (examples of PDPs are ride-sharing, carpool problem, dial-a-ride problem, and vehicle-sharing). Classical VRPs and PDPs can be combined in the so-called \textit{people and freight integrating transportation} problems, which deal with the integration of passenger and freight transportation. Their objective is to increase the occupancy rate by letting the spare seats in the vehicles be used to transport goods (see, e.g., \cite{BBeirigo_2018}, \cite{WChen_2016}, \cite{BLi_2016}, and \cite{VGhilas_2013}). In this way, each vehicle can carry passengers, goods, or both of them.
	
	The \textit{pollution-routing problem}, introduced by \cite{TBekta_2011}, extends the classical VRP by taking into account not only the traveled distance between origin and destination, but also the fuel consumed and the emissions generated by the vehicles. The aim is to find a proper route and speed for each vehicle, allowing customers' requests to be met, by minimizing the overall operational and environmental cost, while respecting time windows and capacity constraints (see, e.g., \cite{RKramer_2014}, \cite{JQian_REglese_2014}, \cite{RFukasawa_QHe_FSantos_YSong_2017}, \cite{MNasri_2018}, and \cite{ISung_PNielsen_2020}). More challenging problems arise when the stochastic (\cite{URitzinger_JPuchinger_RHartl_2016}, \cite{JOyola_HArntzen_DWoodruff_2016a}, \cite{JOyola_HArntzen_DWoodruff_2016b}, \cite{REshtehadi_MFathian_EDemir_2017}), dynamic (\cite{VPillac_MGendreau_CGueret_AMedaglia_2013}, \cite{GBerbeglia_JFCordeau_GLaporte_2010}), and multi-objective (\cite{AGarciaNajera_JBullinaria_2009}, \cite{KGhoseiri_SGhannadpour_2010}, \cite{EDemir_TBektas_GLaporte_2014}, \cite{Kumar_Kondapaneni_Dixit_Goswami_Thakur_Tiwari_2016}) versions are taken into account. 
	When route is considered fixed and the only variable is the vehicle speed, the problem is called \emph{the speed optimization problem} (see \cite{KFagerholt_2010}). The problem considered in our paper can be included in the class of pollution-routing problems because the environmental-impact objective is considered together with speed optimization.
	
	In the literature, few works address a VRP (on a complete graph) by adopting a multi-level formulation, and all of them propose a bilevel problem (see, e.g., \cite{AGupta_2015}, \cite{YMarinakis_AMigdalas_PPardalos_2007}, \cite{YMarinakis_2008}, and \cite{YMa_2014}). In particular, in \cite{AGupta_2015} and \cite{YMarinakis_AMigdalas_PPardalos_2007}, the authors use an assignment-routing formulation to solve a classical~VRP. Similarly, \cite{YMarinakis_2008} propose two nested optimization levels to deal with a VRP integrated with a facility location problem but, in contrast with \cite{AGupta_2015} and \cite{YMarinakis_AMigdalas_PPardalos_2007}, the objective functions involved in each level are conflicting with each other. In our paper, we take advantage of the hierarchical structure at stake by proposing an optimization method that alternates between the assignment-routing problem and the speed optimization problem until a satisfactory solution is returned. Differently from the bilevel approach of \cite{AGupta_2015}, which extends the work of \cite{YMarinakis_AMigdalas_PPardalos_2007} to the bi-objective case by considering efficiency and service quality in each level, we also account for the environmental impact as a third objective.
	
	\subsection{Optimization models for transportation on road networks} \label{subsec:background_operations_research_models_road_network}

When considering road networks and, consequently, non-complete graphs (i.e., graphs that do not contain an edge for every pair of nodes), routing problems face additional challenges because key assumptions are not satisfied~(see, e.g.,~\cite{BFleischmann_1985}, \cite{GCornuejols_1985}, \cite{HBenTicha_2018}, \cite{HBenTicha_TWoensel_2021}). In particular, only a subset of nodes are customers since most of the nodes are associated with cross-roads, which do not have a demand to meet. Therefore, only some nodes need to be visited by the vehicles involved in the problem, which is in contrast with the traditional VRP formulations considering each node in the graph as a customer. Moreover, it may not be possible to find a route that visits nodes and edges only once, thus leading to problems with an empty feasible set. Recently, exact approaches to solve routing problems on a subset of nodes have been proposed in~\cite{RRaeesi_KGZografos_2019} and~\cite{BBoyaci_2021} for VRP, \cite{JRodriguezPereira_2019} for the traveling salesman problem (TSP), and~\cite{HBenTicha_NAbsi_2021} for the shortest path problem. Such papers belong to the stream of works focusing on the so-called Steiner TSP, which was addressed for the first time by~\cite{BFleischmann_1985} and~\cite{GCornuejols_1985}. 

On road networks, arcs may be associated with several attributes (distance, cost, time, etc.), and this implies that the shortest path does not coincide with the cheapest or quickest path. Therefore, an additional approach consists in representing the related non-complete graph by using a multi-graph (\cite{TGaraix_CArtigues_2010}), which allows for multiple arcs between each pair of nodes to account for the best paths associated with each attribute. The computational efficiency of multi-graphs has been questioned in~\cite{ALetchford_SNasiri_AOukil_2014} and reaffirmed in~\cite{HBenTicha_2017}. Instead of resorting to a multi-graph, \cite{KGZografos_KNAndroutsopoulos_2008} and \cite{BHuang_2006} define their VRPs on a road-network, and then propose heuristic procedures that aggregate the attributes associated with each arc. They transform the original graph into a complete one by shortest paths, which is also the approach adopted in our work.

Finally, it is important to mention the class of problems known as {\it arc routing problems} (\cite{ACorberan_2015}), which are based on road networks and, unlike VRPs, the demand is located along the edges of the network. As pointed out in~\cite{HBenTicha_2018}, the main difference between arc routing problems and~VRPs on road networks is that in the former the demand is associated with a large subset of arcs, while in the latter the demand is only on a small subset of nodes. 

	\section{A new trilevel multi-objective formulation architecture for vehicle routing problems with speed optimization} \label{sec:formulation}

\begin{figure}[t!]
		\centering
		\includegraphics[width=5truecm]{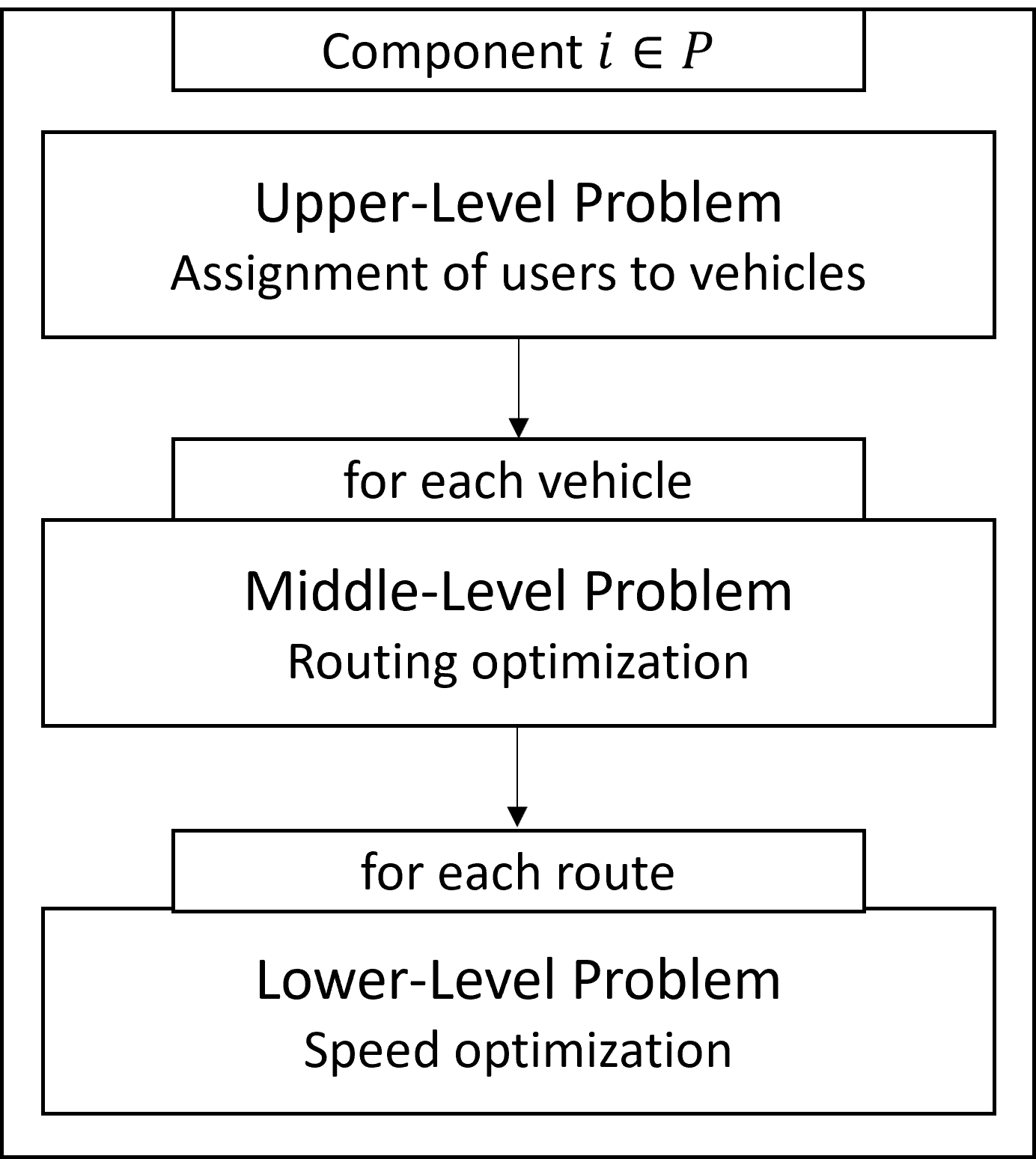}
		\caption{Diagram illustrating the trilevel optimization problem for each transportation component.}
		\label{fig:diag2B}
	\end{figure}

The management of transportation services in smart cities is naturally formulated as an integrated trilevel multi-objective~VRP/speed problem. Figure~\ref{fig:diag2B} presents the main features of the resulting trilevel formulation, which is particularized in Problem~\eqref{prob:general_prob} below. 
	The three hierarchical levels involved in such a problem are repeated for all the $\PP$~possible transportation services. In particular, for each component $i \in \PP$, given a set of vehicles, the Upper-Level (UL) problem determines the optimal assignment of users to such vehicles, represented by a vector~$a^i$ of binary variables. Solving the UL problem requires computing the optimal route associated with each feasible assignment~$a^i$. In fact, given the number of available vehicles and the users assigned to them (this means that $a^i$ is a parameter), the Middle-Level (ML) problem finds an optimal route for each vehicle, represented by a vector~$r^i$ of binary variables. Solving the ML problem requires, in turn, computing the speed (represented by a vector~$v^i$ of real variables) of the given vehicle on each segment of the given route. To this aim, for each possible route ($r^i$~is now handled as a parameter as well), the vehicle speed~$v^i$ is determined by solving the corresponding Lower-Level (LL) problem.
	
	Problem~\eqref{prob:general_prob} depends on several data parameters, among which we highlight~$r^{-i}$ and~$w$. In particular, since the routing decisions~$r^i$ made in each component affect the overall traffic congestion along the streets (and~$v^i$ as a consequence), we introduce the following parameters 
	\begin{equation}\label{eq:notation_j}
	    r^{-i} \; := \; (r^j~|~j \in \{1,\ldots,\PP\} \,\, \text{ and }~j\ne i)
	\end{equation}
	to address the interplay among the different~$\PP$ components, i.e., how routing decisions made in one problem affect the settings in others. The uncertain parameter~$w$ accounts for uncertain factors like the weather or unforeseen events. Furthermore, $w$ can account for the randomness in traffic congestion caused by users that decide not to share their GPS data for privacy concerns (therefore, they are not included in any component).
	\begin{equation}\label{prob:general_prob}
	\setlength{\jot}{5pt}
	\begin{split}
	\min_{\{a^i, r^i, v^i \, | \, i \in \{1,\ldots,\PP\}\}} ~~ & F^i(a^i,r^i,v^i;r^{-i},w) \\[-1ex]
	\mbox{s.t.}~~ & f^i(a^i,r^i,v^i;r^{-i},w)\leq0\\[1.5ex]
	&\, r^i \in \argmin_{r^i, v^i} ~~ G^i(a^i,r^i,v^i;r^{-i},w)\\[-1ex]
	& \qquad\qquad \mbox{s.t.} ~~\, g^i(a^i,r^i,v^i;r^{-i},w)\leq0\\[1.5ex]
	& \qquad\qquad\qquad \, v^i \in \argmin_{v^i} ~~ H^i(a^i,r^i,v^i;r^{-i},w)\\[-1ex]
	& \qquad\qquad\qquad\qquad\qquad \mbox{s.t.} ~~ \, h^i(a^i,r^i,v^i;r^{-i},w)\leq0.
	\end{split}
	\end{equation}
	
	Since the task of the integrated model is to accommodate competing goals to achieve societal benefits, all of the three problems have a multi-objective nature. Accordingly, each objective function is a vector-valued function having three scalar-valued functions as components, each related to one of the three categories of objectives:
	(i) generating a low impact on the environment; (ii) increasing efficiency in providing the transportation services; (iii) meeting requests of users by providing high quality. Functions~$f^i, g^i$, and~$h^i$ represent the constraints. 
	
	The advantage of considering first such a complex trilevel formulation is to comprehensively understand the overall problem complexity and guide the development of subproblems or relaxations and of efficient solution methods. The three levels can be considered separately in a hierarchical order, choosing the most suitable optimization method for each one. Alternatively, the UL and ML~problems can be combined into a single one, leading to the following bilevel problem 
	\begin{equation}\label{prob:general_bilivel_prob_UL}
	\setlength{\jot}{5pt}
	\hspace*{-2cm} \begin{split}
	\min_{\{a^i, r^i, v^i \, | \, i \in \{1,\ldots,\PP\}\}} ~~ & U^i(a^i,r^i,v^i;r^{-i},w) \\[-1ex]
	\mbox{s.t.}~~ & u^i(a^i,r^i,v^i;r^{-i},w)\leq0
	\end{split}
	\end{equation}
	\begin{equation}\label{prob:general_bilivel_prob_LL}
	\hspace*{-2cm} \begin{split}
	\qquad\qquad\qquad\qquad\qquad\qquad\qquad\ \ & \, v^i \in \argmin_{v^i} ~~ L^i(a^i,r^i,v^i;r^{-i},w)\\[-1ex]
	& \qquad\qquad \mbox{s.t.} ~~\, \ell\,^i(a^i,r^i,v^i;r^{-i},w) \leq 0,
	\end{split}
	\end{equation}
    where the (upper level) problem~\eqref{prob:general_bilivel_prob_UL} is an assignment-routing problem and the (lower level) problem~\eqref{prob:general_bilivel_prob_LL} is a speed optimization problem.

	\section{Integration among different transportation components on a road network} \label{sec:mathematical_programming}

	We now describe the approach proposed to model the integration among the different transportation components. In the urban transportation literature, \textit{volume-delay functions}~(VDFs) are commonly used to describe the fundamental relationships between average speed (km/hour) or travel time (hour) on the one hand and traffic flow (vehicles/hour) or density (vehicles/km) on the other hand~(\cite{RKucharski_ADrabicki_2021,JPaszkowski_etal_2021}). One of the most widely adopted~VDFs is the~BPR function proposed by the~\cite{bureau_public_roads_1964}. To the best of our knowledge, the work by~\cite{AVMugayskikh_VVZakharov_TTuovinen_2018} is the only one using the~BPR function in a~VRP problem (without speed optimization). While~\cite{AVMugayskikh_VVZakharov_TTuovinen_2018} has used such a formula to compute the travel time as a function of the traffic flow,~\cite{RKucharski_ADrabicki_2021} have derived an expression to compute the approximate average speed~$v$ on a street as a function of the traffic density. In particular, such an expression is given by
	\begin{equation}\label{eq:VDF_BPR_function}
        v \; = \; \frac{v_0}{1 + \gamma~ (k/k^{\max})^{\eta}},
    \end{equation}
	where~$v_0$ is the free-flow speed (when there is no traffic congestion), $k$ is the traffic density
(number of vehicles per unit distance), $k^{\max}$ is the traffic density when the street is at full capacity, and~$\gamma$ and~$\eta$ are positive parameters. To keep notation simple, we will use the same~$k$ and~$k^{\max}$ to denote numbers of vehicles instead of densities. 
	
	We now want to propose an approximate reformulation of~\eqref{eq:VDF_BPR_function} to relate the average speed on each edge of a non-complete graph to the routing decisions made in other problem components. The non-complete graph associated with a road network is denoted by~$\G=(\N, \E)$, where~$\N$ is the set of nodes representing all the locations relevant for making decisions (which includes not only customer nodes, but also road intersections and connections), while~$\E$ is the set of directed edges connecting pairs of nodes in~$\N$. For the sake of simplicity, we will omit the argument~$i$ associated with the component under consideration, and we will use the  argument~$j$ to denote the parameters that represent the routing decisions made in the other problem components, according to the notation used in~\eqref{eq:notation_j}. Let~${r}_{q\,n_1\,n_2}$ be a binary variable equal to $1$ if vehicle $q \in \V$ traverses~$(n_1,n_2) \in \E$. For each component~$j \in \{1,\ldots,\PP\}$, with~$j \ne i$, we denote the corresponding routing parameters as~$r^j_{q\,n_1\,n_2}$. Let~$k_{n_1\,n_2}$ be the number of vehicles traversing edge~$(n_1\,n_2) \in \E$ that are not controlled by component~$i$, i.e.,
	\begin{equation}\label{eq:traffic_density}
        k_{n_1\,n_2} \; = \; \sum_{j=1, \,j\ne i}^{P}\sum_{q \in {\V}^j} r^j_{q\,n_1\,n_2} + \omega_{n_1\,n_2}.
    \end{equation}
	In~\eqref{eq:traffic_density}, the first term is given by the sum of the routing parameters associated with the other problem components. The second term, represented by~$\omega_{n_1 \, n_2}$, is a non-negative integer random parameter representing the number of vehicles traversing the edge~$(n_1,n_2)$ that are associated with users not sharing their GPS data for privacy concerns (therefore, they are not included in any component). Moreover, we define~$k^{\max}_{n_1\,n_2}$ as the maximum number of vehicles that can traverse edge~$(n_1,n_2)$ at the same time, which can be estimated dividing~$d_{n_1\,n_2}$ by the average length of the vehicles in the smart city streets. 
	
	Denoting the varying speed bounds on edge~$(n_1,n_2)$ as~$v_{n_1 \, n_2}^{\min}$ and $v_{n_1 \, n_2}^{\max}$, introduced for the sake of clarity,
    the speed upper bound is affected by the traffic congestion according to the relation
    \begin{equation}\label{eq:integration}
        v_{n_1\,n_2}^{\max} \; = \; \max\left\{v_{\,n_1\,n_2}^{\maxlim} \left( 1 + \gamma \left(\frac{k_{n_1\,n_2}}{k^{\max}_{n_1\,n_2}}\right)^{\eta}\right)^{-1}, \, v_{\,n_1\,n_2}^{\minlim}\right\},
    \end{equation}
which is inspired by~\eqref{eq:VDF_BPR_function}. The speed lower bound~$v_{n_1\,n_2}^{\min}$ is set equal to~$v_{\,n_1\,n_2}^{\minlim}$ for all~$(n_1, n_2) \in \E$. To ensure~$k_{n_1\,n_2} \le k^{\max}_{n_1\,n_2}$, one should also include the following constraint 
\begin{equation}
\sum_{q \in \V} r_{q\,n_1\,n_2} + k_{n_1\,n_2} \; \le \; k^{\max}_{n_1\, n_2}, \quad \forall (n_1, n_2) \in {\E}. \label{ML_constr_capacity_edge}
\end{equation} 
In~\eqref{eq:integration}, the maximum speed limit~$v_{\,n_1\,n_2}^{\maxlim}$ plays the same role as the free-flow speed in~\eqref{eq:VDF_BPR_function}. In particular, such a speed limit is decreased based on the ratio between~$k_{n_1\,n_2}$ and~$k^{\max}_{n_1\,n_2}$. 
Therefore, the more vehicles traverse an edge, the slower the speed on that edge is. The~$\max$ function is used to ensure that the speed upper bound~$v_{n_1\,n_2}^{\max}$ is greater than the minimum speed limit~$v_{\,n_1\,n_2}^{\minlim}$. Since in practice vehicles may traverse a given edge at different times, the maximum speed resulting from formula~\eqref{eq:integration} gives a pessimistic upper bound on the actual speed of vehicles. 
    

	\section{Formulating and solving an integrated VRP/speed model for freight transportation on a road network}\label{sec:VRP_complete_graph}

The transportation services (also called problem components) arising in a smart city can be categorized into~VRPs, PDPs, or a combination of the two (\textit{People and Freight Integrating Transportation Problems} (PFITP)). Although~PDPs can be considered variants of the classical~VRP, we denote as~VRPs all the problem components where a set of customers wait at fixed locations for the delivery of orders, while with~PDPs we refer to the components characterized by people or objects that need to be picked up at their origins and dropped off at their destinations. In our paper, we provide an integrated formulation for a~VRP/speed problem concerning a specific VRP-type component, i.e., freight transportation, which is sufficient to illustrate all our contributions.  

\begin{figure}
    \centering
    \includegraphics[scale=0.50]{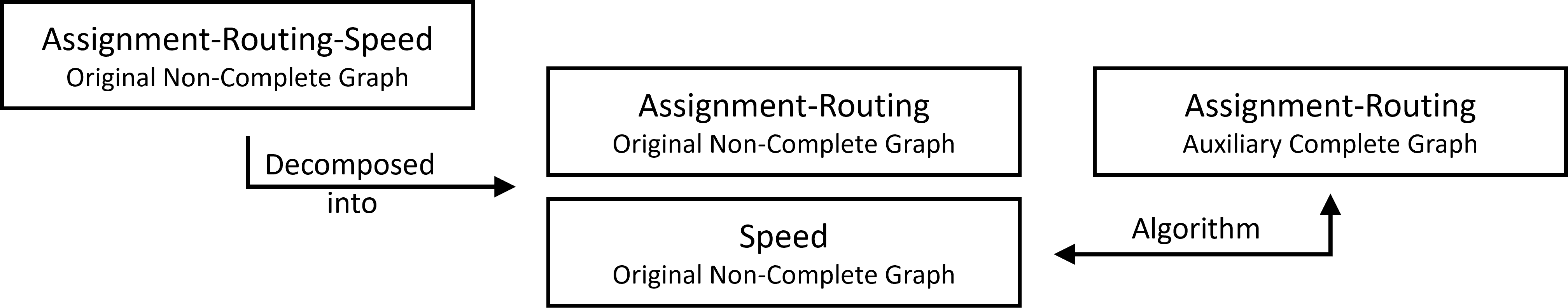}
    \caption{Diagram of the approach proposed to solve a VRP/speed problem on a road network (non-complete graph).}
    \label{fig:greenRoadsDiagram}
\end{figure}

After introducing the notation used to formulate an integrated (bilevel)~VRP/speed model for freight transportation on a road network (Subsection~\ref{sec:notation_2}), we describe the environmental, efficiency, and service objective functions (Subsection~\ref{sec:objective_functions}). Then, we present the assignment-routing problem (Subsection~\ref{sec:extending_assignrouting_road_networks}) and the speed optimization problem (Subsection~\ref{sec:speed_optimization}) that are obtained by decomposing a~VRP/speed problem on a non-complete graph according to upper and lower level problems in~\eqref{prob:general_bilivel_prob_UL}--\eqref{prob:general_bilivel_prob_LL}. Finally, to solve the~VRP/speed problem considered, we propose an optimization algorithm that alternates between the assignment-routing problem on an auxiliary complete graph and the speed optimization problem on the original non-complete graph (Subsection~\ref{sec:algorithm}). The resulting approach is illustrated in Figure~\ref{fig:greenRoadsDiagram}.

	\subsection{Notation for an integrated VRP/speed model for freight transportation}\label{sec:notation_2}

 \begin{table}
    \small
    \centering
    \begin{tabular}{ l|l } 
     \hline
    $\N$ & Set of nodes representing all the locations relevant for making decisions. \\ 
     $\N_{\sub}$ & Subset of $\N$ composed of customer nodes (${\N}_{\sub} \subseteq \N$). \\
     $\E$ & Set of directed edges connecting pairs of nodes in~$\N$.  \\ 
     ${\E}_{\sub}$ & Set of directed edges connecting each pair of customer nodes in~${\N}_{\sub}$, i.e., \\
     & $\{(n_1,n_2) \in {\N}_{\sub} \times {\N}_{\sub} ~|~ n_1 \ne n_2\}.$\\
     $\V$ & Set of available vehicles, i.e.,~$\{1,\ldots,\veh\}$, with $\veh$ maximum number of vehicles. \\
     ${\N}^{q}(r)$ & Set of nodes visited by vehicle~$q \in \V$ based on the routing vector~$r$, i.e.,\\
     & $\{n \in \N~|~r_{q \, n \, \bar{n}} = 1 \text{ for some } (n,\bar{n}) \in \E\}$.\\[2pt]
      $\N^{q}_{\sub}(a)$ & Set of customer nodes visited by vehicle~$q \in \V$ based on the assignment vector~$a$, i.e.,\\
     & $\{n \in {\N}_{\sub}~|~a_{q \, n} = 1\}$.\\
      ${\E}^{q}(r)$ & Set of edges traversed by vehicle~$q \in \V$ based on the routing vector~$r$, i.e.,\\
     & $\{(n_1,n_2) \in \E~|~r_{q \, n_1 \, n_2} = 1\}$.\\[2pt]
      $\E^{q}_{\ordy}(r)$ & Ordered set of the edges in~$\E$ (possibly repeated) traversed by vehicle~$q \in \V$ \\ & based on the routing vector~$r$.\\[2pt]
      ${\E}_{\sub}^{q}(r^{\C})$ & Set of edges in~${\E}_{\sub}$ traversed by vehicle~$q \in \V$ based on a routing vector~$r^{\C}$\\
     & defined on a complete graph, i.e., $\{(n_1,n_2) \in {\E}_{\sub}~|~r^{\C}_{q \, n_1 \, n_2} = 1\}$.\\
     $\SP(\bar{n}_1,\bar{n}_2)$ & Sequence of edges in~$\E$ corresponding to the shortest path on the non-complete \\ 
     & graph between the customer nodes~$\bar{n}_1 \in {\N}_{\sub}$ and~$\bar{n}_2 \in {\N}_{\sub}$.\\[2pt]
     \hline
     $\PP$ & Number of transportation services (problem components). \\
     $d_{n_1\,n_2}$ & Length of the edge~$(n_1,n_2) \in \E$.  \\[2pt]
     $v^{\maxlim}_{n_1\,n_2}$ & Maximum speed limit on the edge~$(n_1,n_2) \in \E$.\\[2pt]
     $v^{\minlim}_{n_1\,n_2}$ & Minimum speed limit on the edge~$(n_1,n_2) \in \E$.\\
     dep & Node in $\N$ corresponding to the depot. \\
     cop & Node in~$\N$ corresponding to a copy of the depot (it can be visited more than once). \\
     $f_q(v)$ & Speed-dependent function computing the emissions (per unit distance) \\ 
     & generated by vehicle $q \in \V$ when traveling at speed~$v$.\\
     $r^j_{q\,n_1\,n_2}$ & Binary parameter equal to~1 if vehicle $q$ from component~$j$ \\ & traverses edge~$(n_1,n_2) \in \E$, with $j \in \{1,\ldots,\PP\}$ and $j \ne i$.\\
     $\omega_{n_1 \, n_2}$ & Positive integer random parameter representing the number of vehicles traversing \\ & the edge $(n_1,n_2) \in \E$ that are not included in any component.\\
     $k_{n_1\,n_2}$ & Number of vehicles traversing the edge~$(n_1,n_2) \in \E$ not controlled by component~$i$.\\
     $k^{\max}_{n_1\,n_2}$ & Maximum number of vehicles that can traverse edge~$(n_1,n_2) \in \E$ at the same time.\\[2pt]
     $v_{n_1 \, n_2}^{\max}$ & Speed upper bound on the edge~$(n_1,n_2) \in \E$.\\
     $v_{n_1 \, n_2}^{\min}$ & Speed lower bound on the edge~$(n_1,n_2) \in \E$.\\[2pt]
     $d^{\SP}_{n_1\,n_2}$ & Length of the shortest path~$\SP(\bar{n}_1,\bar{n}_2)$ for all~$(\bar{n}_1,\bar{n}_2) \in {\E}_{\sub}$. \\
     $\C_q$ & Capacity of vehicle~$q \in \V$. \\
     $\SC_q$ & Setup cost of vehicle~$q \in \V$. \\
     $\DEM_n$ &  Demand of customer $n \in \N_{\sub}$. \\
     $[a_n,b_n]$ & Time window for delivery at customer~$n \in \N_{\sub}$. \\
     $p_n$ & Penalty cost incurred when customer~$n \in \N_{\sub}$ is visited out of the time window. \\[2pt]
     \hline
     $c_1^q(p,r)$ & Function that returns the~$p$-th element of~$\E^{q}_{\ordy}(r)$, with~$p \in \{1,\ldots,|{\E}^{q}_{\ordy}(r)|\}$.\\[2pt]
     $c_2^q(n,p,r)$ & Function that reads~$\E^{q}_{\ordy}(r)$ to count the number of times node~$n \in \N$ \\ & has been visited by vehicle~$q \in \V$ until the edge associated with the~$p$-th \\ & element is traversed.\\[2pt]
     $c_3^q(n_1,n_2,\bar{n}_1,\bar{n}_2)$ & Function that reads the shortest path~$\SP(\bar{n}_1,\bar{n}_2)$ to count the number of times \\ & edge~$(n_1,n_2) \in \E$ is visited by vehicle~$q \in \V$ in~$\SP(\bar{n}_1,\bar{n}_2)$. \\[2pt]
     \hline
    \end{tabular}
    \caption{List of sets, parameters, and functions used in the integrated~VRP/speed model for freight transportation.}\label{tab:list_var_1}
    \end{table}
    
\begin{table}
    \small
    \centering
    \begin{tabular}{ l|l } 
     \hline
     $a_{q\,n}$ & Binary variable equal to~$1$ if vehicle $q \in \V$ is assigned to customer $n \in \N_{\sub}$. \\ 
     ${r}_{q\,n_1\,n_2}$ & Binary variable equal to $1$ if vehicle $q \in \V$ traverses~$(n_1,n_2) \in \E$.  \\ 
     $r^{\rep}_{q\,n_1\,n_2}$ & Integer variable representing the number of times that vehicle $q \in \V$ \\ & traverses~$(n_1,n_2) \in \E$ after the first time. \\
     $v_{q\,n_1\,n_2}$ & Speed of vehicle~$q \in \V$ on the edge $(n_1,n_2) \in \E$. \\
     $t_{q\,n_1\,n_2}$ & Time spent by vehicle~$q \in \V$ to traverse $(n_1,n_2) \in \E$.\\
     $\emis_{q\,n_1\,n_2}$ & Emissions generated by vehicle~$q \in \V$ to traverse $(n_1,n_2) \in \E$.\\
     ${\arr}_{q\,n}^{i}$ & Non-negative real variable representing the arrival time of vehicle~$q \in \V$ at \\ & node~$n \in \N$ when it is visited for the $i$-th time, with~$i = c_2^q(n,p,r)$, $\forall q \in \V$, \\ & $n \in {\N}^{q}(r)$, and $p \in \{1,\ldots,|{\E}^{q}_{\ordy}(r)|\}$.\\[2pt]
     \hline
    \end{tabular}
    \caption{List of optimization variables used in the integrated~VRP/speed model for freight transportation.}\label{tab:list_param}
    \end{table}

    Tables~\ref{tab:list_var_1}--\ref{tab:list_param} introduce the sets, parameters, functions, and optimization variables required to formulate an integrated~VRP/speed model for freight transportation. Throughout this section, the following notation is used to denote the vectors of assignment, routing, and speed variables: 
    \begin{equation}\label{eq:notation_arv}
        \begin{alignedat}{2}
            &a \; := \; (a_{q \, n} \; | \; q \in {\V}, \, n \in {\N}_{\sub}),\\
            &r \; := \; \left((r_{q \, n_1 \, n_2} \; | \; q \in {\V}, \, (n_1, n_2) \in \E), \, (r^{\rep}_{q \, n_1 \, n_2} \; | \; q \in {\V}, \, (n_1, n_2) \in \E)\right), \\
            &v \; := \; \left((v_{q \, n_1 \, n_2} \; | \; q \in {\V}, \, (n_1, n_2) \in \E), \, ({\arr}^i_{q\,n} \; | \; q \in {\V}, \, n \in \N )\right),
        \end{alignedat}
    \end{equation}
    where~$(\textbf{u},\textbf{v})$ is adopted to denote the concatenation of two vectors~$\textbf{u}$ and~$\textbf{v}$. We will use the superscript~$\C$ when the vectors defined in~\eqref{eq:notation_arv} refer to a problem defined on a complete graph.

	We point out that each customer node in the set~$\N_{\sub} \subseteq \N$ is associated with a demand~$\DEM_n$ for a product (for simplicity, we assume single commodity), with~$n \in \N_{\sub}$. Only the nodes in~$\N_{\sub}$ need to be visited by the vehicles involved in the problem, as opposed to the traditional VRP formulations considering each node in the graph as a customer node (see Subsection~\ref{subsec:background_operations_research_models_road_network}).	
	\subsection{Environmental, efficiency, and service objective functions}\label{sec:objective_functions}
Objective function~\eqref{obj_funct_environment} below represents the environmental impact, which is given by the total generated emissions~($\TGE$). Objective function~\eqref{obj_funct_efficiency} below represents the efficiency, which is given by the total setup cost~($\TSC$), total driving time~($\TDT$), and maximum driving distance~($\MDD$). To compute~$\TSC$, we have included a copy of the depot in the set~$\N$ (denoted as~$\cop$) and the edge~$(\dep,\cop)$ in the set~$\E$, and we have replaced all the edges~$(\dep,n)$ in~$\E$ with~$(\cop,n)$ to ensure that a vehicle visits the copy after leaving the depot and before visiting other nodes. This ensures that the setup cost is only paid once for each vehicle regardless of the number of times the copy of the depot is visited. To compute~$\MDD$, we consider the total number of times each edge~$(n_1,n_2)$ in~$\E$ is traversed, which is given by the sum of~$r_{q\,n_1\,n_2}$ and~$r^{\rep}_{q\,n_1\,n_2}$. Objective function~\eqref{obj_funct_quality} below represents the service quality, which is given by the late arrival cost~($\LAC$) and maximum arrival time~($\MAT$), where the non-negative real variable~${\arr}^{1}_{q\,n}$ represents the arrival time of vehicle~$q \in \V$ at node $n \in \N$ when it is visited for the first time. We point out that in the computational experiments the terms in objective functions~\eqref{obj_funct_efficiency}--\eqref{obj_funct_quality} are normalized to avoid numerical issues due to their different units of measurement.
	\begin{alignat}{3}
	& \text{Environmental impact: } && \quad \TGE \label{obj_funct_environment}\\ \vspace{1ex}
	& \text{Efficiency: } && \quad \TSC + \TDT + \MDD \label{obj_funct_efficiency}\\ \vspace{1ex}
    & \text{Service quality: } && \quad \LAC + \MAT \label{obj_funct_quality}
    \end{alignat}
    \vspace{-0.4cm}
    \begin{alignat*}{2}
    & \TGE = \sum_{q \in \V} \sum_{\substack{(n_1, n_2) \in \E}} {\emis}_{q\,n_1\,n_2} 
    && \TSC = \sum_{q \in {\V}} {\SC}_q \, r_{q \dep \cop}\\
    &\TDT = \sum_{q \in {\V}} \sum_{\substack{(n_1, n_2) \in \E}} t_{q\,n_1\,n_2}
    && \MDD = \max_{\substack{q \in V}} \left\{ \sum_{(n_1,n_2) \in \E} (r_{q\,n_1\,n_2} + r^{\rep}_{q\,n_1\,n_2}) \, d_{n_1\,n_2}\right\}\\
    &\LAC = \sum_{q \in {\V}} \sum_{n\in {\N}_{\sub}} \max \{ {\arr}^{1}_{q\,n}-b_n, 0\} \, p_n \qquad
    && \MAT = \max_{\substack{n \in N_{\sub} \\ q \in V}} \{{\arr}^{1}_{q\,n}\} 
    \end{alignat*}

Note that the max functions used in the objective functions~\eqref{obj_funct_efficiency}--\eqref{obj_funct_quality} can be linearized using additional non-negative variables~(i.e.,~$z$, $u$, $w$), re-writing~$\MDD$, $\LAC$, and~$\MAT$ as 
	\begin{equation}\label{eq:obj_lin}
    {\MDD}_{\ell in} = w, \qquad\quad {\LAC}_{\ell in} = \sum_{q \in {\V}} \sum_{n\in {\N}_{\sub}} z_{q\,n} \, p_n, \qquad\quad {\MAT}_{\ell in} = u, 
    \end{equation}
	and considering the following constraints.
	\begin{alignat}{3}
	    & w && \ge \sum_{(n_1,n_2) \in \E} (r_{q\,n_1\,n_2} + r^{\rep}_{q\,n_1\,n_2}) \, d_{q\,n_1\,n_2}, \quad &&\forall q \in \mbox{V} \label{LL_constr_linearization_route}\\
	    & z_{q \, n} && \ge {\arr}^{1}_{q\,n} - b_n, \quad && \forall q \in \mbox{V}, \, n \in N_{\sub} \label{LL_constr_linearization_late_arrivals}\\
	    & u && \ge {\arr}^{1}_{q\,n}, \quad &&\forall q \in \mbox{V}, \, n \in N_{\sub} \label{LL_constr_linearization_arrivals}
	\end{alignat}
	For the rest of the paper, in similar situations where max functions appear, one can linearize them as it was done in~\eqref{eq:obj_lin}--\eqref{LL_constr_linearization_arrivals}.

\subsection{The assignment-routing problem}\label{sec:extending_assignrouting_road_networks}	
When considering road networks, which are represented by non-complete graphs, the standard~VRP/speed formulations based on complete graphs can easily lead to infeasible problems since feasible assignment-routing-speed solutions may not exist. In particular, each customer node may be visited by more than one vehicle and/or more than once by the same vehicle. Moreover, a vehicle may need to both visit some nodes (including the depot) and traverse some edges more than once.

    In the assignment-routing problem, we minimize objective functions~\eqref{obj_funct_environment}--\eqref{obj_funct_quality} to determine an assignment-routing solution~$(a,r)$ on a non-complete graph~$\G(\N,\E)$. The resulting vector of routing variables~$r$ is associated with a sequence of edges in~$\E$ that starts and ends at the depot and possibly visits nodes and edges more than once. According to the notation introduced in Table~\ref{tab:list_var_1}, let~$\E^{q}_{\ordy}(r)\subseteq\E$ be an ordered set of the edges (possibly repeated) traversed by vehicle~$q$, where edges are sorted based on the order in which they are traversed along the route given by~$r$. We can now introduce two functions based on~$\E^{q}_{\ordy}(r)$. Function~$c_1^q(p,r)$ returns the~$p$-th element of~$\E^{q}_{\ordy}(r)$, with~$p \in \{1,\ldots,|{\E}^{q}_{\ordy}(r)|\}$, while function~$c_2^q(n,p,r)$ reads~$\E^{q}_{\ordy}(r)$ to count the number of times node~$n$ has been visited by vehicle~$q$ until the edge associated with the~$p$-th element is traversed.
	Moreover, we define the optimization variable~${\arr}_{q\,n}^{i}$ as the arrival time of vehicle~$q$ at node~$n$ when it is visited for the $i$-th time, i.e.,
	\begin{equation}\label{eq:variable_arr}
	{\arr}_{q\,n}^{i} \ge 0, \text{ with } i = c_2^q(n,p,r), \quad \forall q \in \V, \,n \in {\N}^{q}(r), \, p \in \{1,\ldots,|{\E}^{q}_{\ordy}(r)|\}.
	\end{equation}
	By considering the speed and arrival time variables as parameters, we can write the assignment-routing problem on the non-complete graph~$\G({\N},{\E})$ as a particular case of problem~\eqref{prob:general_bilivel_prob_UL} for our choice of freight transportation component (we use~``(NC)'' to denote that we are here considering a non-complete graph):
	
	{\small
	\begin{equation}\label{prob:assignment-routing-prob-NC}
	\begin{alignedat}{4}
	    & \text{Assignment-Routing (NC):} \quad  & & \text{$U_1 = \TGE$, $U_2 = \TSC + \TDT + \MDD$, and~$U_3 = \LAC + \MAT$.} \\
	    & & & \text{Constraint function~$u$ is given by~\eqref{UL_constr_assign_1_mod}--\eqref{variables_NC} below.}
	\end{alignedat}
	\end{equation}}

Constraint~\eqref{UL_constr_assign_1_mod} below ensures that each customer node is assigned to at least one vehicle (not necessarily just one). Constraint~\eqref{UL_constr_cap_1_mod} below ensures that the demand assigned to a vehicle does not exceed the vehicle capacity. Constraints~\eqref{UL_constr_assign_2_mod} and~\eqref{UL_constr_assign_3_mod} below ensure that each customer node is visited at least once by the assigned vehicle. Constraint~\eqref{ML_constr_flow_start_mod} below ensures that all the vehicles in~$\V$ are used and their routes start at the depot. Constraint~\eqref{ML_constr_flow_conserv_1_mod} below requires flow conservation at each node. Constraint~\eqref{ML_constr_capacity_edge_NC} below is analogous to~\eqref{ML_constr_capacity_edge} and accounts for integration with other components.  Constraints~\eqref{LL_constr_time_1_assignrouting} and~\eqref{LL_constr_emis_1_assignrouting} below compute the travel time and generated emissions for each edge, considering the total number of times the edge is traversed. In particular, $d_{n_1n_2}/v_{q\,n_1\,n_2}$ represents the time required to travel along the edge~$(n_1,n_2)$, while $f_q(v_{q\,n_1\,n_2}) \, d_{n_1 \, n_2}$ computes the emissions generated along the same edge. Constraint~\eqref{eq:precedence_assignrouting} below extends the classical precedence constraint to the case with nodes visited more than once.  
	    \begin{alignat}{1}
	    & \mathbox{A}{\sum_{q \in {\V}}} a_{q\,n} \; \ge \; 1, \quad \forall n \in {\N}_{\sub} \label{UL_constr_assign_1_mod}\\\vspace{2ex}
    	& \mathbox{A}{\sum_{n \in {\N}_{\sub}}} a_{q\,n} {\DEM}_{n} \leq {\C}_q, \quad \forall q \in {\V} \label{UL_constr_cap_1_mod}\\\vspace{2ex}
    	& \mathbox{A}{\sum_{(n_1,n_2) \in \E}} (r_{q\,n_1\,n_2} + r^{\rep}_{q\,n_1\,n_2})\; \ge \; a_{q\,n_1}, \quad \forall q \in {\V}, n_1 \in {\N}_{\sub} \label{UL_constr_assign_2_mod}\\\vspace{2ex}
	    & \mathbox{A}{\sum_{(n_1,n_2) \in \E}} (r_{q\,n_1\,n_2} + r^{\rep}_{q\,n_1\,n_2}) \; \ge \; a_{q\,n_2}, \quad \forall q \in {\V}, n_2 \in {\N}_{\sub} \label{UL_constr_assign_3_mod} \\\vspace{2ex}
    	& \quad r_{q\,\dep\,\cop} \; = \; 1, \quad \forall q \in {\V} \label{ML_constr_flow_start_mod}\\\vspace{2ex} 
    	& \mathbox{A}{\sum_{(n_1,n_2) \in {\E}_{\sub}}} (r_{q\,n_1\,n_2} + r^{\rep}_{q\,n_1\,n_2}) - \sum_{(n_2,n_1) \in {\E}_{\sub}} (r_{q\,n_2\,n_1} + r^{\rep}_{q\,n_2\,n_1}) = 0, \quad \forall q \in {\V},\, n_2 \in {\N}_{\sub} \label{ML_constr_flow_conserv_1_mod}\\\vspace{2ex}
    	& \mathbox{A}{\sum_{q \in \V}} r_{q\,n_1\,n_2} + k_{n_1\,n_2} \; \le \; k^{\max}_{n_1\, n_2}, \quad \forall (n_1, n_2) \in {\E} \label{ML_constr_capacity_edge_NC} \\\vspace{2ex}
    	& \quad t_{q\,n_1\,n_2} = (r_{q\,n_1\,n_2} + r^{\rep}_{q\,n_1\,n_2}) \, d_{n_1\,n_2}/v_{q\,n_1\,n_2}, \quad \forall q \in {\V}, \, (n_1, n_2) \in {\E} \label{LL_constr_time_1_assignrouting} \\\vspace{2ex}
    	& \quad {\emis}_{q\,n_1\,n_2} = (r_{q\,n_1\,n_2} + r^{\rep}_{q\,n_1\,n_2}) \, d_{n_1\,n_2} \, f_q(v_{q\,n_1\,n_2}), \quad  \forall q \in {\V}, \, (n_1, n_2) \in {\E} \label{LL_constr_emis_1_assignrouting}  \\\vspace{2ex}
    	& \quad {\arr}_{q\,n_1}^{i} - {\arr}_{q\,n_2}^{j} + t_{q\,n_1\,n_2} \le M(1-r_{q\,n_1\,n_2}), \nonumber\\
    	& \qquad\qquad\qquad\qquad\qquad\qquad \text{ with } \; (n_1,n_2) = c_1^q(p,r), \; i = c_2^q(n_1,p,r), \; j = c_2^q(n_2,p,r) \nonumber\\
    	& \qquad\qquad\qquad\qquad\qquad\qquad~ \forall q \in \V \text{ and } p \in \{1,\ldots,|{\E}^{q}_{\ordy}(r)|\}\label{eq:precedence_assignrouting}\\\vspace{2ex}
    	& \quad a_{q\,\bar{n}} \in \{0,1\}, \, r_{q\,n_1\,n_2} \in \{0,1\}, \, r^{\rep}_{q\,n_1\,n_2} \in \mathbb{Z}^{+}, \, t_{q\,n_1\,n_2} \in \mathbb{R}, \, {\emis}_{q\,n_1\,n_2} \in \mathbb{R}, \forall q \in \mbox{V}, \nonumber\\\vspace{2ex}
    	& \quad  \bar{n} \in {\N}_{\sub}, \, (n_1, n_2) \in {\E}, \text{ and } {\arr}^{i}_{q\,n_1}, \, {\arr}^{j}_{q\,n_2} \ge 0 \text{ in \eqref{eq:precedence_assignrouting} } \, \forall q \in \mbox{V}, \, p \in \{1,\ldots,|{\E}^{q}_{\ordy}(r)|\} \label{variables_NC}
	    \end{alignat}

\begin{figure}
    \centering
        \subfloat[]{
        \begin{tikzpicture}[node distance={20mm}, thick, main/.style = {draw, circle},scale=0.33]
            \node[main] (0) {$0$};
            \node[main] (6) [right of=0] {$6$};
            \node[main] (7) [below of=6] {$7$}; 
            \node[main,fill={rgb:black,1;white,9}] (8) [below of=0] {$8$};
            \node[main] (5) [left of=0] {$5$}; 
            \node[main] (2) [left of=5] {$2$};
            \node[main,fill={rgb:black,1;white,9}] (1) [below of=2] {$1$};
            \node[main,fill={rgb:black,1;white,9}] (3) [below of=5] {$3$};
            \draw[<->] (0) -- (6);
            \draw[->] (8) -- (6);
            \draw[<->] (8) -- (7);
            \draw[->] (6) -- (7);
            \draw[<->] (0) -- (5);
            \draw[->] (5) -- (2);
            \draw[<->] (2) -- (1);
            \draw[->] (2) -- (3);
            \draw[<->] (3) -- (5);
        \end{tikzpicture}
        }
        \hspace{2cm}   
        \subfloat[]{
        \begin{tikzpicture}[node distance={15mm}, thick, main/.style = {draw, circle},scale=0.33] 
            \node[main] (0) {$0$};
            \node[main,fill={rgb:black,1;white,9}] (8) [below right of=0] {$8$};
            \node[main,fill={rgb:black,1;white,9}] (1) [below left of=0] {$1$};
            \node[main,fill={rgb:black,1;white,9}] (3) [below right of=1] {$3$};
            \draw[<->] (0) -- (8);
            \draw[<->] (0) -- (1);
            \draw[<->] (0) -- (3);
            \draw[<->] (1) -- (3);
            \draw[<->] (3) -- (8);
            \draw[<->] (1) -- (8);
        \end{tikzpicture}
        }
        \caption{Example showing the approach used to obtain an auxiliary complete graph from a non-complete graph. The non-complete graph $\G({\N},{\E})$ reported in~(a) is transformed into the complete graph $\G({\N}_{\sub},{\E}_{\sub})$ in~(b) solely based on the customer nodes.}
        \label{fig:example}
    \end{figure}
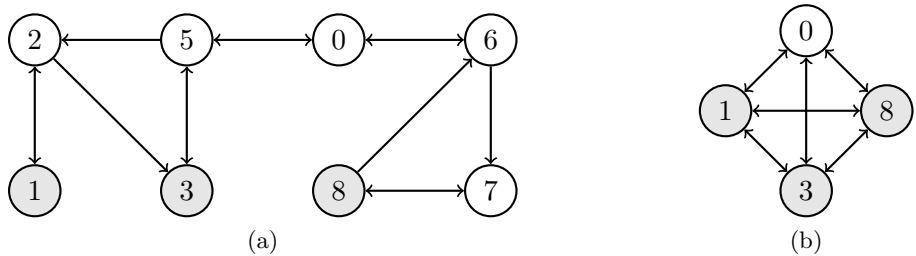

	In our approach, to obtain an approximate assignment-routing solution for problem~\eqref{prob:assignment-routing-prob-NC}, we first solve the assignment-routing problem on an auxiliary complete \linebreak graph~$\G({\N}_{\sub},{\E}_{\sub})$ with set of nodes given by~${\N}_{\sub}$ and set of edges given by~${\E}_{\sub} = \{(n_1,n_2) \in {\N}_{\sub} \times {\N}_{\sub} ~|~ n_1 \ne n_2\}$. In particular, the original road network is transformed into a complete graph by considering the lengths of the shortest paths between each pair of customer nodes as the lengths of the new edges (such lengths can be computed efficiently by using, for example, the Dijkstra algorithm proposed in~\cite{EWDijkstra_1959}). The approach proposed is illustrated in~Figure~\ref{fig:example}. The non-complete graph~$\G(\N, \E)$ in~(a) represents a small road network with single vehicle where the depot is labeled as~0, the customer nodes are denoted as~1, 3, and~8, and all the edges are assumed to have a unitary length except $(2,3)$, whose weight is equal to $4$. To transform such a graph into the complete graph~$\G({\N}_{\sub}, {\E}_{\sub})$ in~(b), we need to compute the shortest paths between the depot and each customer node and between each pair of customer nodes, which are reported in Table~\ref{tab:example} along with their lengths. Note that the~distances between the start and end nodes of the shortest paths are used as lengths for the edges in the complete graph. Although this approach can be computationally expensive when the number of nodes in the original graph is large as compared to the number of edges, as is the case for urban road networks, we pursue this strategy since in our case only a relatively small subset of nodes are customers and, therefore, the resulting complete graph can be computed relatively quickly. 
	
	\begin{table}[]
	    \centering
	    \begin{tabular}{c|c|c}
	    $\bar{n}_1$, $\bar{n}_2$ &  $\SP(\bar{n}_1, \bar{n}_2)$ & $d^{\SP}_{n_1\,n_2}$\\ \hline
	         0, 8 & $\{(0, 6), (6, 7), (7, 8)\}$ & 3\\
	         8, 0 & $\{(8, 6), (6, 0)\}$ & 2\\
	         0, 3 & $\{(0, 5), (5, 3)\}$ & 2\\
	         3, 0 & $\{(3, 5), (5, 0)\}$ & 2\\
	         0, 1 & $\{(0, 5), (5, 2), (2, 1)\}$ & 3\\
	         1, 0 & $\{(1, 2), (2, 3), (3, 5), (5, 0)\}$ & 7\\
	         8, 3 & $\{(8, 6), (6, 0), (0, 5), (5, 3)\}$ & 4\\
	         3, 8 & $\{(3, 5), (5, 0), (0, 6), (6, 7), (7, 8)\}$ & 5\\
	         3, 1 & $\{(3, 5), (5, 2), (2, 1)\}$ & 3\\
	         1, 3 & $\{(1, 2), (2, 3)\}$ & 5\\
	         1, 8 & $\{(1, 2), (2, 3), (3, 5), (5, 0), (0, 6), (6, 7), (7, 8)\}$ & 10\\
	         8, 1 & $\{(8, 6), (6, 0), (0, 5), (5, 2), (2, 1)\}$ & 5
	    \end{tabular}
        \caption{Shortest path and corresponding length between the depot and each customer node and between each pair of customer nodes in the non-complete graph~$\G({\N},{\E})$ considered in Figure~\ref{fig:example}.}
	    \label{tab:example}
	\end{table}
	
	After building the auxiliary complete graph~$\G({\N}_{\sub},{\E}_{\sub})$, standard~VRP techniques can be used to obtain an optimal assignment-routing on such a graph by solving problem~\eqref{prob:assignment-routing-prob} below (``(C)'' denotes that we are considering a complete graph), where the optimization variables are the ones listed in Table~\ref{tab:list_var_auxiliary}. In such a problem, we assume that each vehicle can visit the nodes (including the depot) and edges at most once.

\begin{table}
    \small
    \centering
    \begin{tabular}{ l|l } 
     \hline
     $a^{\C}_{q\,n}$ & Binary variable equal to~$1$ if vehicle $q \in \V$ is assigned to customer $n \in \N_{\sub}$. \\[2pt] 
     $r^{\C}_{q\,n_1\,n_2}$ & Binary variable equal to $1$ if vehicle $q \in \V$ traverses edge~$(n_1,n_2) \in {\E}_{\sub}$. \\[2pt]  
     $v^{\C}_{q\,n_1\,n_2}$ & Speed of vehicle~$q \in \V$ on the edge $(n_1,n_2) \in {\E}_{\sub}$. \\[2pt] 
     $\arr^{\C}_{q\,n}$ & Arrival time of vehicle~$q \in \V$ at node $n \in {\N}_{\sub}$. \\[2pt] 
     $t^{\C}_{q\,n_1\,n_2}$ & Time spent by vehicle~$q \in \V$ to traverse $(n_1,n_2) \in {\E}_{\sub}$. \\[2pt] 
     $\emis^{\C}_{q\,n_1\,n_2}$ & Emissions generated by vehicle~$q \in \V$ to traverse $(n_1,n_2) \in {\E}_{\sub}$. \\[2pt]
     \hline
    \end{tabular}
    \caption{List of the optimization variables used in the formulation of the assignment-routing problem on the auxiliary complete graph.}\label{tab:list_var_auxiliary}
    \end{table}
	
	{\small
	\begin{equation}\label{prob:assignment-routing-prob}
	\begin{alignedat}{4}
	    & \text{Assignment-Routing (C):} \quad  & & \text{$U^{\C}_1 = \overline{\TGE}$, $U^{\C}_2 = \overline{\TSC} + \overline{\TDT} + \overline{\MDD}$, and~$U^{\C}_3 = \overline{\LAC} + \overline{\MAT}$.} \\
	    & & & \text{Constraint function~$u^{\C}$ is given by~\eqref{UL_constr_assign_1}--\eqref{variables} below.}  
	\end{alignedat}
	\end{equation}}
	\begin{alignat*}{2}
    & \overline{\TGE} = \sum_{q \in \V} \sum_{\substack{(n_1, n_2) \in {\E}_{\sub}}} {\emis}^{\C}_{q\,n_1\,n_2} 
    && \overline{\TSC} = \sum_{q \in {\V}} {\SC}_q \sum_{(dep,n) \in {\E}_{\sub}} r^{\C}_{q\,dep\,n}\\
    &\overline{\TDT} = \sum_{q \in {\V}} \sum_{\substack{(n_1, n_2) \in {\E}_{\sub}}} t^{\C}_{q\,n_1\,n_2}
    && \overline{\MDD} = \max_{\substack{q \in V}} \left\{ \sum_{(n_1,n_2) \in {\E}_{\sub}} r^{\C}_{q\,n_1\,n_2} \, d^{\SP}_{n_1\,n_2}\right\}\\
    &\overline{\LAC} = \sum_{q \in {\V}} \sum_{n\in {\N}_{\sub}} \max \{ {\arr}^{\C}_{q\,n}-b_n, 0\} \, p_n \qquad
    && \overline{\MAT} = \max_{\substack{n \in N_{\sub} \\ q \in V}} \{{\arr}^{\C}_{q\,n}\} 
    \end{alignat*}
	
	Constraints \eqref{UL_constr_assign_1}-\eqref{ML_constr_flow_conserv_1} below are standard in the VRP literature and represent the user assignment (each customer is served by one vehicle), capacity constraint (demand cannot exceed vehicle capacity), route assignment (each customer is included in a route), starting flow (each vehicle departs from the depot), and flow conservation (incoming flow is equal to outgoing flow).
	Constraints~\eqref{LL_constr_time_1_C} and~\eqref{LL_constr_emis_1_C} below compute the travel time and the generated emissions on each edge. In \eqref{LL_constr_prec_1_C} below, $M$ is a sufficiently large positive constant, whose value can be set equal to the sum of the largest values of ${\arr}^{\C}_{q\,n_1}$ and $t^{\C}_{q\,n_1\,n_2}$. 
    \begin{align}
	& \mathbox{A}{\sum_{q \in {\V}}} a^{\C}_{q\,n} = 1, \quad \forall n \in {\N}_{\sub} \label{UL_constr_assign_1}\\\vspace{2ex}
	& \mathbox{A}{\sum_{n \in {\N}_{\sub}}} a^{\C}_{q\,n} {\DEM}_{n} \leq {\C}_q, \quad \forall q \in {\V} \label{UL_constr_cap_1}\\\vspace{2ex}
	&  \mathbox{A}{\sum_{(n_1,n_2) \in {\E}_{\sub}}} r^{\C}_{q\,n_1\,n_2} = a^{\C}_{q\,n_1}, \quad \forall q \in {\V}, \, n_1 \in {\N}_{\sub} \label{UL_constr_assign_2} \\\vspace{2ex}
	& \mathbox{A}{\sum_{(n_1,n_2) \in {\E}_{\sub}}} r^{\C}_{q\,n_1\,n_2} = a^{\C}_{q\,n_2}, \quad \forall q \in {\V}, \,n_2 \in {\N}_{\sub} \label{UL_constr_assign_3}\\\vspace{2ex}
	& \mathbox{A}{\sum_{(dep,n) \in {\E}_{\sub}}} r^{\C}_{q\,dep\,n} = 1, \quad \forall q \in {\V} \label{ML_constr_flow_start}\\\vspace{2ex}
	& \mathbox{A}{\sum_{(n_1,n_2) \in {\E}_{\sub}}} r^{\C}_{q\,n_1\,n_2} - \sum_{(n_2,n_1) \in {\E}_{\sub}}r^{\C}_{q\,n_2\,n_1} = 0, \quad \forall q \in {\V},\, n_2 \in {\N}_{\sub} \label{ML_constr_flow_conserv_1}\\\vspace{2ex}
	& \quad t^{\C}_{q\,n_1\,n_2} = (r^{\C}_{q\,n_1\,n_2} \, d^{\SP}_{n_1\,n_2})/v^{\C}_{q\,n_1\,n_2}, \quad \forall q \in {\V}, \, (n_1, n_2) \in {\E}_{\sub} \label{LL_constr_time_1_C} \\\vspace{4ex}
	& \quad {\emis}^{\C}_{q\,n_1\,n_2} = r^{\C}_{q\,n_1\,n_2} \, d^{\SP}_{n_1\,n_2} \, f_q(v^{\C}_{q\,n_1\,n_2}), \quad  \forall q \in {\V}, \, (n_1, n_2) \in {\E}_{\sub} \label{LL_constr_emis_1_C}  \\\vspace{4ex}
    & \quad {\arr}^{\C}_{q\,n_1} - {\arr}^{\C}_{q\,n_2} + t^{\C}_{q\,n_1\,n_2} \le M(1-r^{\C}_{q\,n_1\,n_2}), 
	\quad \forall q \in {\V}, \,(n_1, n_2) \in {\E}_{\sub}, \, n_2 \ne \dep \label{LL_constr_prec_1_C} \\\vspace{4ex}
	& \quad a^{\C}_{q\,n} \in \{0,1\}, \, r^{\C}_{q\,n_1\,n_2} \in \{0,1\}, \, {\arr}^{\C}_{q\,n} \ge 0, \, t^{\C}_{q\,n_1\,n_2} \in \mathbb{R}, \, {\emis}_{q\,n_1\,n_2} \in \mathbb{R}, \nonumber\\\vspace{4ex}
	& \quad \forall q \in \mbox{V},\, n \in {\N}_{\sub}, \, (n_1, n_2) \in {\E}_{\sub} \label{variables}
	\end{align}
	Note that the assignment-routing problem~\eqref{prob:assignment-routing-prob} is solely composed of binary variables (the variables~$t^{\C}_{q\,n_1\,n_2}$ and~$\emis^{\C}_{q\,n_1\,n_2}$ are only introduced for the sake of clarity in the presentation of the formulation). 

After solving the auxiliary problem~$\eqref{prob:assignment-routing-prob}$, an approximate solution~$(\tilde{a},\tilde{r})$ for the original assignment-routing problem~\eqref{prob:assignment-routing-prob-NC} can be derived from the shortest paths associated with each edge in the optimal routing on the auxiliary complete graph. In the example illustrated in Figure~\ref{fig:example}, assume that an optimal route for the VRP defined on the auxiliary complete graph~$\G({\N}_{\sub},{\E}_{\sub})$ is~0-8-3-1-0. Then, a route for the original non-complete graph~$\G({\N},{\E})$ can be obtained by considering the shortest paths between~0-8, 8-3, 3-1, and~1-0, which leads to 0-6-7-8-6-0-5-3-5-2-1-2-3-5-0. 

More formally, according to Table~\ref{tab:list_var_1}, let us define the shortest path on the non-complete graph between two customer nodes as a sequence of edges in~$\E$ denoted by~$\SP(\bar{n}_1,\bar{n}_2)$, with~$(\bar{n}_1, \bar{n}_2) \in {\E}_{\sub}$.
The function~$c_3^q(n_1,n_2,\bar{n}_1,\bar{n}_2)$ reads~$\SP(\bar{n}_1,\bar{n}_2)$ to count the number of times edge~$(n_1,n_2) \in \E$ is visited by vehicle~$q \in \V$. Moreover, given an optimal solution~$(\bar{a}^{\C},\bar{r}^{\C})$ to the assignment-routing problem~\eqref{prob:assignment-routing-prob}, let~${\E}^{q}_{\sub}(\bar{r}^{\C}) \subseteq {\E}_{\sub}$ be the set of edges on the auxiliary complete graph traversed by vehicle~$q$, with~$q \in \V$, based on the routing vector~$\bar{r}^{\C}$, i.e.,
	\[
	{\E}^{q}_{\sub}(\bar{r}^{\C}) ~=~ \{(n_1,n_2) \in {\E}_{\sub}~|~\bar{r}^{\C}_{q \, n_1 \, n_2} = 1\}. 
	\] 
For all~$q \in \V$, $n \in {\N}_{\sub}$, and $(n_1,n_2) \in \E$, we can set
\begin{equation}\label{eq:solution_recovery}
\begin{alignedat}{3}
&\tilde{a}_{q\,n} \; &&= \; \bar{a}^{\C}_{q\,n}, \\
&\tilde{r}_{q\,n_1\,n_2} \; &&= \; 
\begin{cases}
1 \quad \text{ if } (n_1,n_2) \in \SP(\bar{n}_1,\bar{n}_2) \text{ for some } (\bar{n}_1, \bar{n}_2) \in {\E}_{\sub}^q(\bar{r}^{\C}), \\
0 \quad \text{ otherwise}, 
\end{cases} \\
&\tilde{r}^{\rep}_{q\,n_1\,n_2} \; &&= \; 
\left(\sum_{(\bar{n}_1,\bar{n}_2) \in {\E}_{\sub}^q(\bar{r}^{\C}) } c^q_3(n_1,n_2,\bar{n}_1,\bar{n}_2)\right)-1.
\end{alignedat}
\end{equation}

Note that constraints~\eqref{UL_constr_assign_1}--\eqref{ML_constr_flow_conserv_1} in the auxiliary assignment-routing problem~\eqref{prob:assignment-routing-prob} ensure that constraints~\eqref{UL_constr_assign_1_mod}--\eqref{ML_constr_flow_conserv_1_mod} in the original assignment-routing problem~\eqref{prob:assignment-routing-prob-NC} are satisfied at the point~$(\tilde{a}, \tilde{r})$ resulting from~\eqref{eq:solution_recovery}. In the formulation of problem~\eqref{prob:assignment-routing-prob}, we do not include a constraint enforcing~\eqref{ML_constr_capacity_edge_NC} because it would not allow solving problem~\eqref{prob:assignment-routing-prob} by using existing~VRP methods, which are not designed to handle such a constraint. Therefore, 
to ensure that constraint~\eqref{ML_constr_capacity_edge_NC} evaluated at~$(\tilde{a},\tilde{r})$ is satisfied for all~$q \in \V$ and~$(n_1, n_2) \in \E$, we use a penalization approach. In particular, when such a constraint is violated for some vehicles~$\hat{q} \in \V$ and edges~$(\hat{n}_1, \hat{n}_2) \in \E$, we solve the auxiliary assignment-routing~\eqref{prob:assignment-routing-prob} again by introducing penalty terms in the objective function to ensure that the edges~$(\hat{n}_1, \hat{n}_2)$ will not be included in the vehicle's route, i.e.,
\begin{equation}\label{eq:penalization}
\Phi(r^{\C}_{\hat{q} \, \bar{n}_1 \, \bar{n}_2}) \; = \; \begin{cases}
 +\infty &\text{ if } r^{\C}_{\hat{q} \, \bar{n}_1 \, \bar{n}_2} = 1,\\[1ex]
 0 &\text{ if } r^{\C}_{\hat{q} \, \bar{n}_1 \, \bar{n}_2} = 0,
\end{cases}
\end{equation}
where~$(\bar{n}_1, \bar{n}_2) \in {\E}_{\sub}$ is such that~$(\hat{n}_1, \hat{n}_2) \in \SP(\bar{n}_1, \bar{n}_2)$.

\subsection{The speed optimization problem}\label{sec:speed_optimization}

	When solving the speed optimization problem, the assignment and routing variables on the non-complete graph, represented by the vectors~$a$ and~$r$, respectively, are considered parameters.  
	In particular, according to the notation introduced in Table~\ref{tab:list_var_1}, let~$\N^{q}(r)\subseteq\N$ be the set of nodes visited by vehicle~$q$, with~$q \in \V$, based on the routing vector~$r$, i.e., 
	\begin{alignat*}{3}
	&{\N}^{q}(r) &&~=~ \{n \in \N~|~r_{q \, n \, \bar{n}} = 1 \text{ for some } (n,\bar{n}) \in \E\}. 
	\end{alignat*}
	Similarly, let~$\N^{q}_{\sub}(a)\subseteq\N$ be the set of customer nodes visited by vehicle~$q$, with~$q \in \V$, based on the assignment vector~$a$, i.e., 
	\begin{alignat*}{3}
	&{\N}_{\sub}^{q}(a) &&~=~ \{n \in {\N}_{\sub}~|~a_{q \, n} = 1\}. 
	\end{alignat*}
	Finally, let~$\E^{q}(r)\subseteq\E$ be the set of edges traversed by vehicle~$q$, with~$q \in \V$, based on the routing vector~$r$, i.e.,
	\[
	{\E}^{q}(r) ~=~ \{(n_1,n_2) \in \E~|~r_{q \, n_1 \, n_2} = 1\}. 
	\] 
	
	We can now write the speed optimization problem on the non-complete graph~$\G({\N},{\E})$ as a particular case of problem~\eqref{prob:general_bilivel_prob_LL} for our choice of freight transportation component
	\begin{equation}\label{prob:speed-prob-NC}
	\begin{alignedat}{4}
	    & \text{Speed (NC):} \qquad\qquad& & \text{$L_1 = \TGE$, $L_2 = \TDT$, and~$L_3 = \LAC + \MAT$.} \\
	    & & & \text{$\ell$ is given by~\eqref{LL_constr_time_1}--\eqref{variables_2} below}.
    \end{alignedat}
    \end{equation}
    The terms in the objective functions of problem~\eqref{prob:speed-prob-NC} are defined in~\eqref{obj_funct_environment}--\eqref{obj_funct_quality} (note that~$\E$ can be replaced by~${\E}^q(r)$ and~$ {\N}_{\sub}$ can be replaced by~${\N}_{\sub}^q(a)$). Constraints~\eqref{LL_constr_time_1}--\eqref{eq:precedence} below are analogous to~\eqref{LL_constr_time_1_assignrouting}--\eqref{eq:precedence_assignrouting}. Constraints~\eqref{LL_constr_speed_bound_1} and~\eqref{LL_constr_speed_bound_2} below restrict the value of the vehicle speed on each edge, where~$v^{\max}_{n_1\,n_2}$ is given by~\eqref{eq:integration}, while~$v_{n_1\,n_2}^{\min}$ is set equal to~$v_{\,n_1\,n_2}^{\minlim}$ for all~$(n_1, n_2) \in {\E}^q(r)$.
    \begin{align}
	& \quad t_{q\,n_1\,n_2} = (r_{q\,n_1\,n_2} + r^{\rep}_{q\,n_1\,n_2}) \, d_{n_1\,n_2}/v_{q\,n_1\,n_2}, \quad \forall q \in {\V}, \, (n_1, n_2) \in {\E}^q(r) \label{LL_constr_time_1} \\\vspace{2ex}
	& \quad {\emis}_{q\,n_1\,n_2} = (r_{q\,n_1\,n_2} + r^{\rep}_{q\,n_1\,n_2}) \, d_{n_1\,n_2} \, f_q(v_{q\,n_1\,n_2}), \quad  \forall q \in {\V}, \, (n_1, n_2) \in {\E}^q(r) \label{LL_constr_emis_1}  \\\vspace{2ex}
	& \quad {\arr}_{q\,n_1}^{i} - {\arr}_{q\,n_2}^{j} + t_{q\,n_1\,n_2} \le M(1-r_{q\,n_1\,n_2}), \nonumber\\
	& \qquad\qquad\qquad\qquad\qquad\qquad \text{ with } \; (n_1,n_2) = c_1^q(p,r), \; i = c_2^q(n_1,p,r), \; j = c_2^q(n_2,p,r) \nonumber\\
	& \qquad\qquad\qquad\qquad\qquad\qquad~ \forall q \in \V \text{ and } p \in \{1,\ldots,|{\E}^{q}_{\ordy}(r)|\}\label{eq:precedence}\\\vspace{2ex}
	& \quad v_{q\,n_1\,n_2} \leq v_{n_1\,n_2}^{\max}, \quad \forall q \in {\V}, \, (n_1, n_2) \in {\E}^q(r)  \label{LL_constr_speed_bound_1}\\\vspace{1ex}
	& \quad v_{q\,n_1\,n_2} \geq v_{n_1\,n_2}^{\min}, \quad \, \forall q \in {\V}, \, (n_1, n_2) \in {\E}^q(r) \label{LL_constr_speed_bound_2}\\\vspace{1ex}
	& \quad v_{q\,n_1\,n_2} \ge 0, \, t_{q\,n_1\,n_2} \in \mathbb{R}, \, {\emis}_{q\,n_1\,n_2} \in \mathbb{R}, \forall q \in \mbox{V}, \, (n_1, n_2) \in {\E}, n \in {\N}^{q}(r) \nonumber\\\vspace{2ex}
	& \qquad\quad \text{ and } {\arr}^{i}_{q\,n_1}, \, {\arr}^{i}_{q\,n_2} \ge 0 \text{ in \eqref{eq:precedence} } \, \forall q \in \mbox{V}, \, p \in \{1,\ldots,|{\E}^{q}_{\ordy}(r)|\} \label{variables_2}
	\end{align}
	
    Note that all the optimization variables in the speed-optimization problem~\eqref{prob:speed-prob-NC} are real-valued. Given that constraints~\eqref{LL_constr_time_1}--\eqref{LL_constr_emis_1} are nonlinear, such a speed optimization problem is a nonlinear constrained optimization problem.

	\subsection{An algorithm for assignment-routing-speed optimization}
	\label{sec:algorithm}

	To solve the VRP/speed problem on road networks, we propose an optimization method that alternates between the auxiliary assignment-routing problem~\eqref{prob:assignment-routing-prob} and the speed optimization one~\eqref{prob:speed-prob-NC}. To handle the multiple objective functions in such problems, we use the weighted-sum method (\cite{MEhrgott_2005}, \cite{KMiettinen_2012}). Therefore, given non-negative weights~$\alpha_i$, with~$i \in \{1,2,3\}$, such objective functions are weighted into single-objective scalar functions, i.e.,~$U(a,r,v) = \sum_{i=1}^{3} \alpha_i \, U_i(a,r,v)$ and~$L(a,r,v) = \sum_{i=1}^{3} \alpha_i \, L_i(a,r,v)$. 
	
	The schema of the proposed optimization method is reported in Algorithm~\ref{alg:alg_trilevel}, where we use the notation introduced in~\eqref{eq:notation_arv}. To keep track of the values of the variables across the iterations in the two {\it while loops}, we use two subscripts, i.e.,~$a_{k, \, j}$, $r_{k, \, j}$, $v_{k, \, j}$, where $k$ refers to the outer loop and~$j$ to the inner loop. Moreover, the number of vehicles available at iteration~$k$ is denoted by~${\veh}_k$.  
	
	Given an initial number of available vehicles~${\veh}_0$ and an initial speed vector~$v_{0,0}$, the algorithm performs an exhaustive enumeration until either the maximum number of outer iterations~$k_{\max}$ or the maximum number of customer nodes~$|N_{sub}|$ is reached. For each number of vehicles, the algorithm runs an outer cycle until no change occurs in the assignment, routing, and speed variables. At each outer iteration, an assignment-routing on the non-complete graph is obtained transforming the original network into an auxiliary complete graph and determining the optimal assignment-routing on such a graph by using standard~VRP methods. The resulting assignment-routing solution will be passed to the speed optimization problem as a parameter in order to obtain a speed value for each edge traversed by the vehicle on the road network. More specifically, given the current number of available vehicles~${\veh}_k$ and the current speed vector~$v_{k,j}$, Step~2.1 determines~$\bar{a}^{\C}_{k,j}$ and~$\bar{r}^{\C}_{k,j}$ by solving the assignment-routing problem on the auxiliary complete graph~\eqref{prob:assignment-routing-prob}. Step~2.2 determines an approximate solution~$(\tilde{a}_{k,\,j}, \tilde{r}_{k, \, j})$ to the assignment-routing problem~\eqref{prob:assignment-routing-prob-NC} by applying~\eqref{eq:solution_recovery}. Step~2.3 determines~$\tilde{v}_{k,j}$ by solving the speed optimization problem~\eqref{prob:speed-prob-NC} given the assignment and routing determined at the previous step. The outer loop proceeds until either the maximum number of inner iterations is reached or the current solution is equal to the previous one, meaning that no further improvement is possible. The final solution resulting from the application of the algorithm is denoted as~$(\doublehat{a}_{k}, \, \doublehat{r}_{k}, \, \doublehat{v}_{k})$.
	
	 \begin{algorithm}[hbtp]
	\caption{}\label{alg:alg_trilevel}
	\begin{algorithmic}[1]
		\medskip
        \item[] {\bf Input:} Initial number of vehicles ${\veh}_0 = 1$, initial speed $v_{0, \, 0}$, $N_{\sub}$, $k = 0$, $k_{\max}$, $j_{\max}$, \item[] \qquad\quad\ \ $\doublehat{U} := + \infty$.\nonumber
		\smallskip
		\item[] {\bf While $k \le k_{\max}$ and ${\veh}_k \le |N_{\sub}|$ \bf do}
		\item[] \quad\quad {\bf Step 1.} $\hat{U} := + \infty$, $j = 1$. 
		 \item[] \quad\quad {\bf While $j \le j_{\max}$ and $(a_{k, \, j}, \, r_{k, \, j}, \, v_{k, \, j}) \ne (a_{k, \, j-1}, \, r_{k, \, j-1}, \, v_{k, \, j-1})$ \bf do} \vspace{0.1cm}
		 \item[] \quad\quad\quad\quad {\bf Step 2.1.} Obtain~$(\bar{a}^{\C}_{k,\,j}, \bar{r}^{\C}_{k, \, j})$ by solving the assignment-routing problem~\eqref{prob:assignment-routing-prob} on the auxiliary complete graph with fixed~${\veh}_k$ and~$v_{k, \, j}$.
		 \item[] \quad\quad\quad\quad {\bf Step 2.2.} Obtain an approximate solution~$(\tilde{a}_{k,\,j}, \tilde{r}_{k, \, j})$ to the assignment-routing problem~\eqref{prob:assignment-routing-prob-NC} by applying~\eqref{eq:solution_recovery}.
 \nonumber
 \item[] \quad\quad\quad\quad {\bf Step 2.3.} Obtain $\tilde{v}_{k, \, j}$ by solving the speed optimization problem~\eqref{prob:speed-prob-NC} with fixed $\tilde{a}_{k, \, j}$ and $\tilde{r}_{k, \, j}$. Set~$(a_{k, \, j}, \, r_{k, \, j}, \, v_{k, \, j}) = (\tilde{a}_{k,\,j}, \tilde{r}_{k, \, j}, \tilde{v}_{k, \, j})$. \nonumber
 \item[] \quad\quad\quad\quad {\bf Step 2.4.} {\bf If} $U(a_{k, \, j}, \, r_{k, \, j}, \, v_{k, \, j}) < \hat{U}$. \nonumber
 \item[] \quad\quad\quad\quad\quad \phantom{{\bf Step 2.4.}} $\hat{U} := U(a_{k, \, j}, \, r_{k, \, j}, \, v_{k, \, j})$. \nonumber
 \item[] \quad\quad\quad\quad\quad \phantom{{\bf Step 2.4.}} $(\hat{a}_{k}, \, \hat{r}_{k}, \, \hat{v}_{k}) := (a_{k, \, j}, \, r_{k, \, j}, \, v_{k, \, j})$. \nonumber
  \item[] \quad\quad\quad\quad {\bf Step 2.5.} $j = j+1$. \nonumber
        \item[] \quad\quad {\bf End do}
		\item[] \quad\quad {\bf Step 3.} {\bf If} $\hat{U} < \doublehat{U}$. \nonumber
		\item[] \quad\quad\quad \phantom{{\bf Step 2.}} $\doublehat{U} := \hat{U}$. \nonumber
		\item[] \quad\quad\quad \phantom{{\bf Step 2.}} $(\doublehat{a}_{k}, \, \doublehat{r}_{k}, \, \doublehat{v}_{k}) := (\hat{a}_{k}, \, \hat{r}_{k}, \, \hat{v}_{k})$.
		\item[] \quad\quad {\bf Step 4.} ${\veh}_{k+1} = {\veh}_k + 1$ and $k = k+1$. \nonumber
		\item[] {\bf End do}
		\par\bigskip\noindent
    	\end{algorithmic}
    \end{algorithm}

    \section{Computational experiments}\label{sec:computational_experiments}
    
    
    Two instances have been randomly generated from a dataset gathering hourly traffic information in the New York City streets from 2010 to 2013 (with more than 95,500 nodes and 260,850 arcs), which has been created by~\cite{BDonovan_2015}. In particular, we consider a first instance (referred to as~{\it small} in the remainder of this section) consisting of a graph with 269 nodes, 656 edges, and 10 customer nodes, and a second instance (indicated as {\it large}) based on a graph with~1130 nodes, 2703 edges, and 30 customer nodes. The small and large instances have been obtained by including all the nodes with longitude in~$[-73.970,-73.941]$ and latitude in~$[40.791,40.810]$ and all the nodes with longitude in~$[-73.970,-73.938]$ and latitude in~$[40.760,40.840]$, respectively. 
    
    \begin{figure}
        \centering
        \includegraphics[scale=0.6]{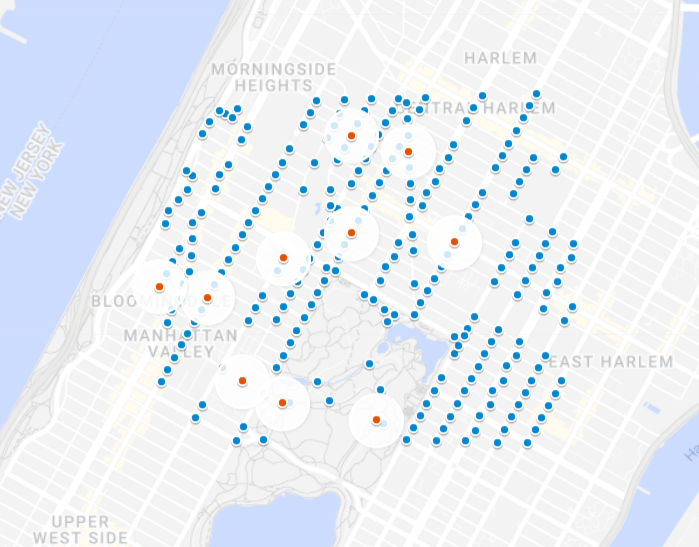}\quad \includegraphics[scale=0.40]{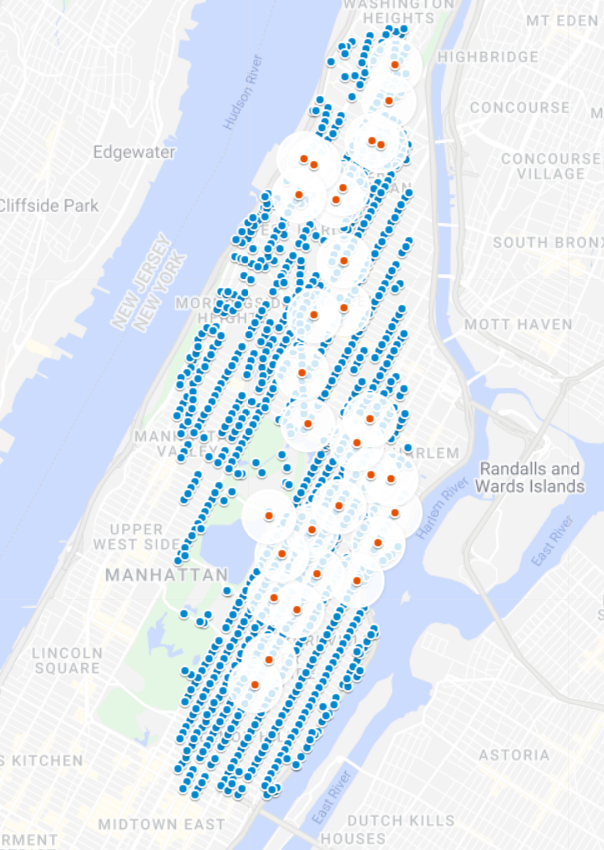}
        \caption{Small (left) and large (right) instances. Customers are represented by the big circles.}
        \label{fig:instances}
    \end{figure}
    
    In all the experiments, the customer nodes in~$N_{\sub}$ have been randomly generated from the graphs. The time window upper bound~$b_n$ for each node~$n \in N_{\sub}$ has been set to $2/3*|N_{\sub}|*\TS(\dep,n)$, where $\TS(\dep,n)$ denotes the time required to reach $n$ from the depot by traveling along the shortest path at a speed of $10$ miles/h (about~$16.09$ km/h), which has resulted in tight windows. Since we are considering an urban setting, the maximum speed limit~$v_{n_1 \, n_2}^{\maxlim}$ has been set to~$25$~miles/h (about~$40.23$ km/h), while the minimum speed limit~$v_{n_1 \, n_2}^{\minlim}$ has been set to~$5$~miles/h (about~$8.05$ km/h). 
    
    The demand parameter~$\DEM_n$ has been set to 0 for all~$n \in \N_{\sub}$ to omit the capacity constraints from these experiments. Both the vehicle setup cost~$\SC_q$ and the cost for late arrivals~$p_n$ have been arbitrarily set to~$\$ 100$ for all~$q \in \V$ and~$n \in \N_{sub}$. The capacity $\capac_{n_1\,n_2}$ of each edge~$(n_1,n_2)$ has been estimated dividing~$d_{n_1 \, n_2}$ by the average length of vehicles in a city, set to~$4.5$~m (which is a reasonable assumption, according to~\cite{VMarkevicius_etal_2017}). Finally, to model the speed-dependent function $f_q$ computing the $\CO_2$ emissions per unit of distance~(gram/km) generated by a vehicle $q \in \V$ traversing an edge~$(n_1,n_2)$, we have used the following formula, reported in~\cite{EDemir_TBektas_GLaporte_review_2014}:
    \begin{equation*}
        f_q(v_{q \, n_1 \, n_2}) = 871 - 16 v_{q \, n_1 \, n_2} + 0.143 v_{q \, n_1 \, n_2}^2 + 32031/v_{q \, n_1 \, n_2}^2,  
    \end{equation*}
    where~$v_{q \, n_1 \, n_2}$ is in km/h.
    
    The integration among the different transportation components is modeled through formula~\eqref{eq:integration} and constraint~\eqref{ML_constr_capacity_edge}, which is analogous to constraint~\eqref{ML_constr_capacity_edge_NC} used in the original assignment-routing problem~\eqref{prob:assignment-routing-prob-NC}. We recall that to ensure that the approximate solution obtained by applying~\eqref{eq:solution_recovery} is feasible with respect to constraint~\eqref{ML_constr_capacity_edge_NC}, we adopt the penalty term~\eqref{eq:penalization}, which requires solving the auxiliary assignment-routing problem~\eqref{prob:assignment-routing-prob} with such a penalty term in the objective function every time such a constraint is violated on the original non-complete graph. 
    According to~\cite{RKucharski_ADrabicki_2021},  $\gamma = 1$ and~$\eta = 2$ are reasonable values for the parameters used in formula~\eqref{eq:integration} when considering urban road networks. Denoting by~$i$ the current component, we handle the variables~$r^j_{q\,n_1\,n_2}$ in~\eqref{eq:integration} and~\eqref{ML_constr_capacity_edge} as parameters, with~$j \in \{1,\ldots,\PP\}$ and~$j \ne i$. 
    In particular, two cases are taken into account to represent the total number of vehicles controlled by other components that traverse an edge~$(n_1, n_2)$. In the first case, which we denote as~{0\%-cap}, no vehicles controlled by other components traverse the edge~$(n_1, n_2)$, while in the second one, denoted as~{50\%-cap}, the number of such vehicles is equal to~50\% of the capacity of the edges (the word \textit{capacity} is used to refer to~$k^{\max}_{n_1\, n_2}$, which is the maximum number of vehicles that can traverse the edge at the same time). Recalling the definition of~$k_{n_1\,n_2}$ in~\eqref{eq:traffic_density}, one has
    \[
    k_{n_1\,n_2} =
    \begin{cases}
     0 & \text{(Case 0\%-cap)}, \\ 0.5 \, k^{\max}_{n_1\, n_2} & \text{(Case 50\%-cap)},
    \end{cases}
    \] 
    for all~$(n_1, n_2) \in \E$.
    
    We recall that in multi-objective optimization one is interested in obtaining a set of points that cannot improve one objective without worsening the values of the other ones. Points with this property are called Pareto optimal solutions (or non-dominated points). To obtain the approximate Pareto fronts, the weighted-sum method with normalization has been applied to all the terms in the objective functions~$U_1$, $U_2$, and $U_3$ of the original assignment-routing problem~\eqref{prob:assignment-routing-prob-NC}, all the terms in~$U_1^{\C}$, $U_2^{\C}$, and $U_3^{\C}$ for the auxiliary assignment-routing problem~\eqref{prob:assignment-routing-prob}, and all the terms in~$L_1$, $L_2$, and~$L_3$ for the speed optimization problem~\eqref{prob:speed-prob-NC}. In particular, given non-negative weights~$\alpha_i$, with~$i \in \{1,\ldots,6\}$, the problems considered in the experiments are~\eqref{prob:assignment-routing-prob-NC-2}--\eqref{prob:speed-prob-NC-2} below. Each term in the objective functions has been normalized using additional weights~$w_i$ or~$\bar{w}_i$, with~$i \in \{1,\ldots,6\}$, which represent the value of each term at the initial solution, i.e.,~$w_1 = \TGE(a_{0,0},r_{0,0},v_{0,0})$, $w_2 = \TSC(a_{0,0},r_{0,0},v_{0,0})$, and similarly for the remaining terms and for the weights~$\bar{w}_i$.
    
    \vspace{-0.1cm}
	{\small
	\begin{equation}\label{prob:assignment-routing-prob-NC-2}
	\begin{alignedat}{4}
	    & \text{Assignment-Routing (NC):} \quad  & & \text{$U = \alpha_1\left(\frac{\TGE}{w_1}\right) + \alpha_2\left(\frac{\TSC}{w_2}\right) + \alpha_3\left(\frac{\TDT}{w_3}\right) +$}\\ & & & \qquad \text{$\alpha_4\left(\frac{\MDD}{w_4}\right) + \alpha_5\left(\frac{\LAC}{w_5}\right) + \alpha_6\left(\frac{\MAT}{w_6}\right)$.} \\
	    & & & \text{Constraint function~$u$ is the same as problem~\eqref{prob:assignment-routing-prob-NC}}.
	\end{alignedat}
	\end{equation}
	\vspace{-0.8cm}
	\begin{equation}\label{prob:assignment-routing-prob-2}
	\begin{alignedat}{4}
	& \phantom{\text{Assignment-Routing (NC):}} \quad  & & \phantom{\text{$U = \alpha_1\left(\frac{\TGE}{w_1}\right) + \alpha_2\left(\frac{\TSC}{w_1}\right) + \alpha_3\left(\frac{\TDT}{w_1}\right) +$}} \\
	& \text{Assignment-Routing (C):} \quad  & & \text{$U^{\C} = \alpha_1\left(\frac{\overline{\TGE}}{\bar{w}_1}\right) + \alpha_2\left(\frac{\overline{\TSC}}{\bar{w}_2}\right) + \alpha_3\left(\frac{\overline{\TDT}}{\bar{w}_3}\right) +$}\\ & & & \qquad \ \ \, \text{$\alpha_4\left(\frac{\overline{\MDD}}{\bar{w}_4}\right) + \alpha_5\left(\frac{\overline{\LAC}}{\bar{w}_5}\right) + \alpha_6\left(\frac{\overline{\MAT}}{\bar{w}_6}\right)$.} \\
	& & & \text{Constraint function~$u^{\C}$ is the same as problem~\eqref{prob:assignment-routing-prob}.} \\
	& \phantom{\text{Assignment-Routing (NC):}} \quad  & & \phantom{\text{$U = \alpha_1\left(\frac{\TGE}{w_1}\right) + \alpha_2\left(\frac{\TSC}{w_1}\right) + \alpha_3\left(\frac{\TDT}{w_1}\right) +$}} \\
    \end{alignedat}
    \end{equation}
    \vspace{-1.8cm}
	\begin{equation}\label{prob:speed-prob-NC-2}
	\begin{alignedat}{4}
	& \phantom{\text{Assignment-Routing (NC):}} \quad  & & \phantom{\text{$U = \alpha_1\left(\frac{\TGE}{w_1}\right) + \alpha_2\left(\frac{\TSC}{w_1}\right) + \alpha_3\left(\frac{\TDT}{w_1}\right) +$}} \\
	& \text{Speed (NC):} & & \text{$L = \alpha_1\left(\frac{\TGE}{w_1}\right) + \alpha_3\left(\frac{\TDT}{w_3}\right) +$}
	\\ & & & \qquad \text{$\alpha_5\left(\frac{\LAC}{w_5}\right) + \alpha_6\left(\frac{\MAT}{w_6}\right)$.} \\
	& & & \text{Constraint function~$\ell$ is the same as problem~\eqref{prob:speed-prob-NC}.} \\
	& \phantom{\text{Assignment-Routing (NC):}} \quad  & & \phantom{\text{$U = \alpha_1\left(\frac{\TGE}{w_1}\right) + \alpha_2\left(\frac{\TSC}{w_1}\right) + \alpha_3\left(\frac{\TDT}{w_1}\right) +$}} \\
    \end{alignedat}
    \end{equation}}
    \vspace{-0.8cm}
    
    Among the solutions returned by the algorithm, we eliminated the dominated points to obtain an approximation of the Pareto front. We point out that although all the terms are considered in the experiments, in the plots we only report the comparison of the approximate Pareto fronts among the conflicting terms in the environmental impact and efficiency objective functions. The terms in the service quality objective function are omitted from the figures since they have been observed not to be conflicting with the other ones. In particular, the maximum driving distance~($\MDD$), which is referred to as {\it Max. Distance} in the figures, has been compared against the total generated~$\CO_2$ emissions~($\TGE$) (referred to as {\it Total Emissions}), the total setup cost~($\TSC$), and total driving time~($\TDT$).
    
    All tests were run on a Linux server with~32GB of RAM and an AMD Opteron 6128 processor running at~2.00 GHz. The assignment-routing problem (Step 2.1 in Algorithm~\ref{alg:alg_trilevel}) has been implemented in~Python~3.7 and solved by using OR-Tools~9.1~(\cite{LPerron_VFurnon_2019}) with default options. To solve the speed optimization problem (Step 2.3), the interior point method implemented in the solver IPOPT has been used~(\cite{AWachter_LTBiegler_2006}) with the parameter~{\it tol}~set to~$10^{-2}$.
 
     \subsection*{Deterministic case}\label{sec:deterministic_instance_without_speed}
    
\begin{figure}
        \centering
        \includegraphics[scale=0.31]{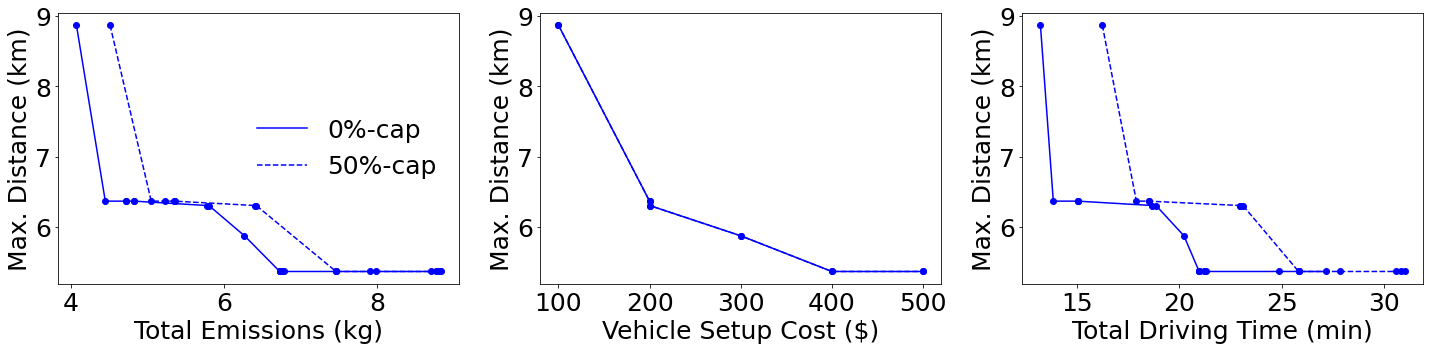}
        \caption{Approximate Pareto fronts for the small instance.}
        \label{fig:general_small}
    \end{figure}

\begin{figure}
        \centering
        \includegraphics[scale=0.31]{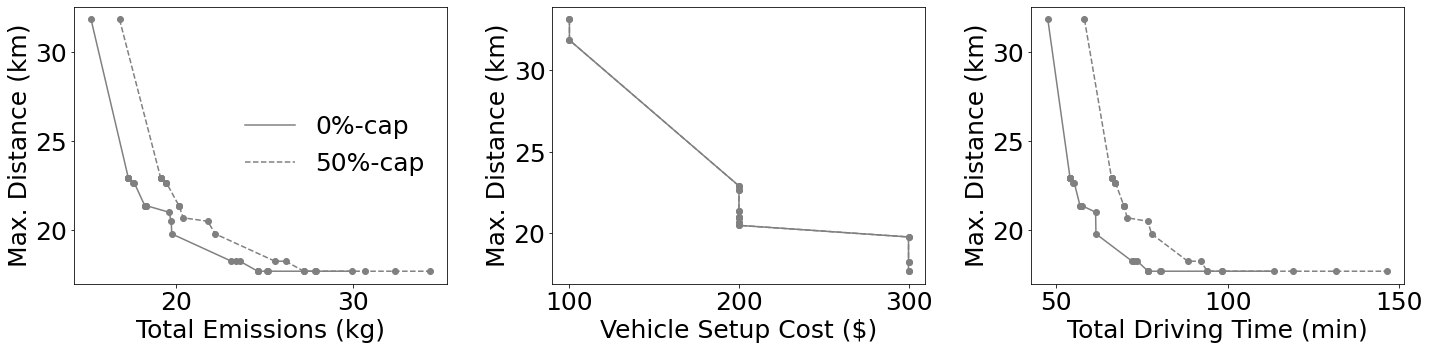}
        \caption{Approximate Pareto fronts for the large instance.}
        \label{fig:general_large}
    \end{figure}

    In the deterministic case, the parameter $\omega_{n_1 \, n_2}$ in~\eqref{eq:integration} has been set to 0 for all~$(n_1,n_2) \in N_{\sub}$. Figures~\ref{fig:general_small}--\ref{fig:general_large} show the comparison between the approximate Pareto fronts obtained for~{0\%-cap} and~{50\%-cap} on the small and large instances. Such figures confirm that the integration among the transportation components has a significant impact on the objectives considered, as illustrated in the plots. Once a set of objective function weights is given, the average CPU time required by Algorithm~\ref{alg:alg_trilevel} on the small instance is~408.58~s and on the large instance is~2995.91~s. The majority of the time is due to the computation of an approximate solution for the speed optimization problem. The development of a heuristic for speed optimization to reduce the computational cost is left for future work.
    
    From the plots related to the small instance, one can observe that reducing the number of vehicles from~5 to~1 (which corresponds to a decrease in the vehicle setup cost from~\$500 to~\$100) leads to an increase of the maximum driving distance from~5.37~km to~8.87~km, which allows achieving~53.27\% (in the 0\%-cap case) and~48.98\% (in the 50\%-cap case) savings in~$\CO_2$ emissions and~51.29\% (0\%-cap) and~47.50\% (50\%-cap) savings in total driving time. Moreover, the first and third plots show that without increasing the maximum driving distance, it is possible to reduce the total emissions from~8.71~kg to~6.73~kg (0\%-cap) and from~8.84~kg to~7.47~kg (50\%-cap) and the total driving time from~27.16~min to~20.97~min (0\%-cap) and from~31.01~min to~25.82~min (50\%-cap). Note that using 5 vehicles does not lead to any significant improvement in terms of maximum driving distance compared to the Pareto solution with 4 vehicles.
    
    On the large instance, only solutions with a number of vehicles between~\$100 and~\$300 are Pareto optimal. Reducing such a number from~3 to~1 (which corresponds to a decrease in the vehicle setup cost from~\$300 to~\$100) leads to an increase of the maximum driving distance from~17.68~km to~31.85~km, and the savings that can be obtained are~49.33\% (0\%-cap) and~51.17\% (50\%-cap) for~$\CO_2$ emissions and~58.17\% (0\%-cap) and~60.27 (50\%-cap) for total driving time. 
    
     \subsection*{Stochastic case}\label{sec:stochastic_instance_with_speed}

    In this subsection, we limit the analysis to the small instance and we consider the~50\%-cap case. The random variable~$\omega_{n_1 \, n_2}$ in~\eqref{eq:integration} has a Poisson distribution with parameter~$\lambda_{n_1 \, n_2}$ for all~$(n_1,n_2) \in N_{\sub}$. Different values of the parameter have been considered based on the capacity of each edge, i.e.,~$\lambda_{n_1 \, n_2} = \beta \, k^{\max}_{n_1 \, n_2}$ with $\beta \in \{0.10,0.15\}$. 
    
\begin{figure}
        \centering
        \includegraphics[scale=0.31]{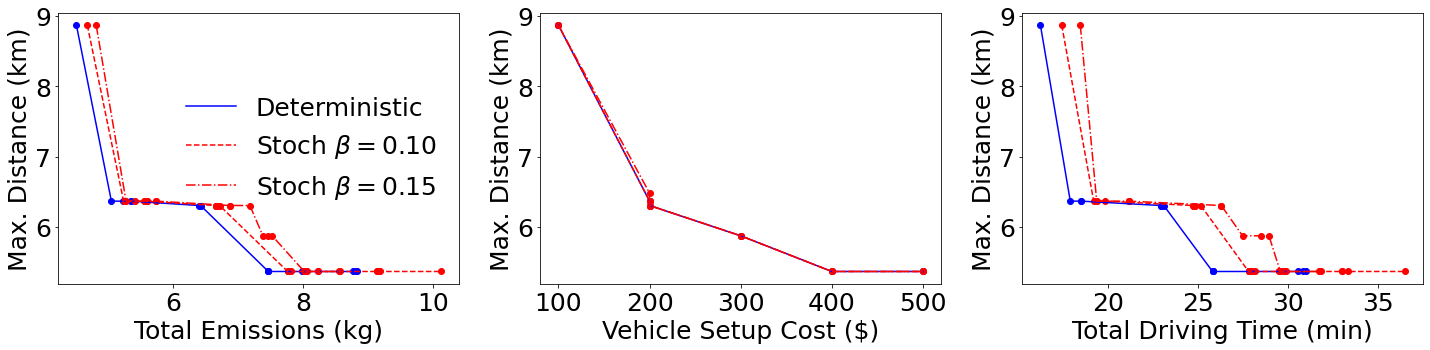}
        \caption{Approximate Pareto fronts for the deterministic case and the stochastic one.}
        \label{fig:stoch}
    \end{figure}

    Figure~\ref{fig:stoch} shows the comparison of the approximate Pareto fronts obtained in the deterministic and stochastic cases. In particular, in the stochastic case, both the~$\CO_2$ emissions and total driving time significantly increase as compared to the deterministic case due to the decrease in speed caused by traffic congestion, according to formula~\eqref{eq:integration} and constraint~\eqref{ML_constr_capacity_edge}. Reducing the number of vehicles from~5 to~1 leads to~53.77\% ($\beta=10\%$) and~43.87\% ($\beta=15\%$) savings in~$\CO_2$ emissions and~52.29\% ($\beta=10\%$) and~42.08\% ($\beta=15\%$) savings in total driving time. The maximum driving distance increases from~5.37~km to~8.87~km. Once a set of objective function weights is given, the average CPU time required by~Algorithm~\ref{alg:alg_trilevel} for~$\beta=10\%$ is~422.54~s and for~$\beta=15\%$ is~637.11~s. 

    \section{Concluding remarks and future work}\label{sec:conclusion}
    
    In this paper, we laid the groundwork for the development of a green-oriented integrated system aimed at managing mobility and transportation services within a smart city. 
    In particular, we proposed a new formulation architecture consisting of three nested optimization levels associated with assignment, routing, and speed decisions. To address the integration among all the different transportation problem components, we developed an innovative constraint on the speed variables (i.e., formula~\eqref{eq:integration} and~constraint~\eqref{ML_constr_capacity_edge}) that allows one to model the impact of the traffic congestion caused by routing decisions made in other components on the component considered (which is a~VRP-type one in our paper). Based on the formulation architecture, we developed formulations for the assignment-routing problem and the speed optimization one for the chosen freight transportation component. The approach proposed to solve the~VRP/speed problem considered computes an approximate solution to the assignment-routing problem on a non-complete graph (i.e., problem~\eqref{prob:assignment-routing-prob-NC}) by first solving such a problem on an auxiliary complete graph (i.e., problem~\eqref{prob:assignment-routing-prob}). The resulting algorithm (i.e.,~Algorithm~\eqref{alg:alg_trilevel}) alternates between the auxiliary assignment-routing problem~\eqref{prob:assignment-routing-prob} and the speed optimization problem~\eqref{prob:speed-prob-NC} on the non-complete graph.
    
	The computational experiments show that the proposed approach is able to determine a set of Pareto optimal solutions among the conflicting terms of the objective functions considered. Moreover, the experimental results show the importance of considering the impact of different components due to traffic congestion. Therefore, the integrated model can be used as a decision support tool to help find the best trade-off among competing criteria.
	
	The development of concrete formulations for~PDP and~PFITP components and the analysis of potential conflicts among the objective functions when such components are considered is left for future work. Moreover, further research is needed to improve the proposed concrete formulation for road networks by removing the assumption that the vehicle speed used on an edge traversed more than once is the same as the first time the edge was traversed. Finally, an additional avenue of research will be focused on the speed optimization problem, which has been so far solved by using an algorithm that provides solutions with local convergence guarantees but, at the same time, high computational cost. Although the use of an exact algorithm provides us with a good solution in terms of vehicle speed given the assignment and routing decisions, in practice the use of heuristic algorithms for speed optimization is needed to ensure a real-time response.   
	
	To manage the transportation network and grant users access to the smart decision-making system, one should develop an application for smart devices. Such an app would acquire data from the users by asking them for the type of service they need (use their own vehicles, car sharing, ride sharing, carpooling, etc.) and other relevant information (such as origin and destination). Moreover, users will need to specify whether to concede the vehicle-related decisions on assignment, routing, and speed to the app or make their own decisions based on private objectives. The acquired data would then be used to feed the integrated optimization model proposed in this paper.


\begin{thebibliography}{65}
	\providecommand{\natexlab}[1]{#1}
	\providecommand{\url}[1]{\texttt{#1}}
	\expandafter\ifx\csname urlstyle\endcsname\relax
	\providecommand{\doi}[1]{doi: #1}\else
	\providecommand{\doi}{doi: \begingroup \urlstyle{rm}\Url}\fi
	
	\bibitem[Acheampong and Cugurullo(2019)]{RAAcheampong_FCugurullo_2019}
	R.~A. Acheampong and F.~Cugurullo.
	\newblock Capturing the behavioural determinants behind the adoption of
	autonomous vehicles: Conceptual frameworks and measurement models to predict
	public transport, sharing and ownership trends of self-driving cars.
	\newblock \emph{Transportation Research Part F: Traffic Psychology and
		Behaviour}, 62:\penalty0 349--375, 2019.
	
	\bibitem[Arroub et~al.(2016)Arroub, Zahi, Sabir, and
	Sadik]{AArroub_BZahi_ESabir_etal_2016}
	A.~Arroub, B.~Zahi, E.~Sabir, and M.~Sadik.
	\newblock A literature review on smart cities: Paradigms, opportunities and
	open problems.
	\newblock In \emph{2016 International Conference on Wireless Networks and
		Mobile Communications (WINCOM)}, pages 180--186, 2016.
	
	\bibitem[Beirigo et~al.(2018)Beirigo, Schulte, and Negenborn]{BBeirigo_2018}
	B.~Beirigo, F.~Schulte, and R.R. Negenborn.
	\newblock Integrating people and freight transportation using shared autonomous
	vehicles with compartments.
	\newblock \emph{IFAC-PapersOnLine}, 51:\penalty0 392--397, 01 2018.
	
	\bibitem[Bektas and Laporte(2011)]{TBekta_2011}
	T.~Bektas and G.~Laporte.
	\newblock The pollution-routing problem.
	\newblock \emph{Transportation Research Part B: Methodological}, 45:\penalty0
	1232--1250, 09 2011.
	
	\bibitem[Ben~Ticha et~al.(2017)Ben~Ticha, Absi, Feillet, and
	Quilliot]{HBenTicha_2017}
	H.~Ben~Ticha, N.~Absi, D.~Feillet, and A.~Quilliot.
	\newblock Empirical analysis for the vrptw with a multigraph representation for
	the road network.
	\newblock \emph{Comput. Oper. Res.}, 88:\penalty0 103–116, dec 2017.
	
	\bibitem[Ben~Ticha et~al.(2018)Ben~Ticha, Absi, Feillet, and
	Quilliot]{HBenTicha_2018}
	H.~Ben~Ticha, N.~Absi, D.~Feillet, and A.~Quilliot.
	\newblock Vehicle routing problems with road-network information: State of the
	art.
	\newblock \emph{Networks}, 72:\penalty0 393--406, 2018.
	
	\bibitem[{Ben Ticha} et~al.(2021{\natexlab{a}}){Ben Ticha}, Absi, Feillet, and
	Quilliot]{HBenTicha_NAbsi_2021}
	H.~{Ben Ticha}, N.~Absi, D.~Feillet, and A.~Quilliot.
	\newblock The steiner bi-objective shortest path problem.
	\newblock \emph{EURO Journal on Computational Optimization}, 9:\penalty0
	100004, 2021{\natexlab{a}}.
	
	\bibitem[{Ben Ticha} et~al.(2021{\natexlab{b}}){Ben Ticha}, Absi, Feillet,
	Quilliot, and Woensel]{HBenTicha_TWoensel_2021}
	H.~{Ben Ticha}, N.~Absi, D.~Feillet, A.~Quilliot, and T.~Woensel.
	\newblock {The Time-Dependent Vehicle Routing Problem with Time Windows and
		Road-Network Information}.
	\newblock \emph{SN Operations Research Forum}, 2:\penalty0 1--25, March
	2021{\natexlab{b}}.
	
	\bibitem[Berbeglia et~al.(2010)Berbeglia, Cordeau, and
	Laporte]{GBerbeglia_JFCordeau_GLaporte_2010}
	G.~Berbeglia, J.F. Cordeau, and G.~Laporte.
	\newblock Dynamic pickup and delivery problems.
	\newblock \emph{European Journal of Operational Research}, 202:\penalty0 8 --
	15, 2010.
	
	\bibitem[Boyacı et~al.(2021)Boyacı, Dang, and Letchford]{BBoyaci_2021}
	B.~Boyacı, T.~Dang, and A.~Letchford.
	\newblock Vehicle routing on road networks: How good is euclidean
	approximation?
	\newblock \emph{Computers \& Operations Research}, 129, 05 2021.
	
	\bibitem[Bronstein(2009)]{ZBronstein_2009}
	Z.~Bronstein.
	\newblock Industry and the smart city.
	\newblock \emph{Dissent}, 56:\penalty0 27--34, 06 2009.
	\newblock \doi{10.1353/dss.0.0062}.
	
	\bibitem[{Bureau of Public Roads}(1964)]{bureau_public_roads_1964}
	{Bureau of Public Roads}.
	\newblock \emph{Traffic Assignment Manual for Application with a Large, High
		Speed Computer}.
	\newblock U. S. Department of Commerce, Bureau of Public Roads, Office of
	Planning, Urban Planning Division, 1964.
	
	\bibitem[Cattaruzza et~al.(2017)Cattaruzza, Absi, Feillet, and
	Gonzalez-Feliu]{DCattaruzza_2017}
	D.~Cattaruzza, N.~Absi, D.~Feillet, and J.~Gonzalez-Feliu.
	\newblock Vehicle routing problems for city logistics.
	\newblock \emph{EURO Journal on Transportation and Logistics}, 6:\penalty0
	51--79, 04 2017.
	
	\bibitem[Chen et~al.(2018)Chen, Mes, and Schutten]{WChen_2016}
	W.~Chen, M.~Mes, and M.~Schutten.
	\newblock Multi-hop driver-parcel matching problem with time windows.
	\newblock \emph{Flexible services and manufacturing journal}, 30:\penalty0
	517--553, 9 2018.
	\newblock Springer deal.
	
	\bibitem[Corber\'{a}n and Laporte(2015)]{ACorberan_2015}
	A.~Corber\'{a}n and G.~Laporte.
	\newblock \emph{Arc Routing: Problems, Methods, and Applications}.
	\newblock Society for Industrial and Applied Mathematics, 2015.
	
	\bibitem[Cornu\'{e}jols et~al.(1985)Cornu\'{e}jols, Fonlupt, and
	Naddef]{GCornuejols_1985}
	G.~Cornu\'{e}jols, J.~Fonlupt, and D.~Naddef.
	\newblock The traveling salesman problem on a graph and some related integer
	polyhedra.
	\newblock \emph{Math. Program.}, 33:\penalty0 1–27, sep 1985.
	
	\bibitem[Cugurullo(2020)]{FCugurullo_2020}
	F.~Cugurullo.
	\newblock Urban artificial intelligence: From automation to autonomy in the
	smart city.
	\newblock \emph{Frontiers in Sustainable Cities}, 2, 2020.
	
	\bibitem[Demir et~al.(2014{\natexlab{a}})Demir, Bektaş, and
	Laporte]{EDemir_TBektas_GLaporte_2014}
	E.~Demir, T.~Bektaş, and G.~Laporte.
	\newblock The bi-objective pollution-routing problem.
	\newblock \emph{European Journal of Operational Research}, 232:\penalty0 464 --
	478, 2014{\natexlab{a}}.
	
	\bibitem[Demir et~al.(2014{\natexlab{b}})Demir, Bektaş, and
	Laporte]{EDemir_TBektas_GLaporte_review_2014}
	E.~Demir, T.~Bektaş, and G.~Laporte.
	\newblock A review of recent research on green road freight transportation.
	\newblock \emph{Eur. J. Oper. Res.}, 237:\penalty0 775--793,
	2014{\natexlab{b}}.
	
	\bibitem[Dijkstra(1959)]{EWDijkstra_1959}
	E.~W. Dijkstra.
	\newblock A note on two problems in connexion with graphs.
	\newblock \emph{Numer. Math.}, 1:\penalty0 269–271, dec 1959.
	
	\bibitem[Donovan(2015)]{BDonovan_2015}
	B.~Donovan.
	\newblock Link level traffic estimates.
	\newblock 2015.
	\newblock URL \url{https://uofi.app.box.com/v/NYC-traffic-estimates}.
	
	\bibitem[Ehrgott(2005)]{MEhrgott_2005}
	M.~Ehrgott.
	\newblock \emph{Multicriteria Optimization}, volume 491.
	\newblock Springer Science \& Business Media, Berlin, 2005.
	
	\bibitem[Eshtehadi et~al.(2017)Eshtehadi, Fathian, and
	Demir]{REshtehadi_MFathian_EDemir_2017}
	R.~Eshtehadi, M.~Fathian, and E.~Demir.
	\newblock Robust solutions to the pollution-routing problem with demand and
	travel time uncertainty.
	\newblock \emph{Transportation Research Part D: Transport and Environment},
	51:\penalty0 351 -- 363, 2017.
	
	\bibitem[Fagerholt et~al.(2010)Fagerholt, Laporte, and
	Norstad]{KFagerholt_2010}
	K.~Fagerholt, G~Laporte, and I.~Norstad.
	\newblock Reducing fuel emissions by optimizing speed on shipping routes.
	\newblock \emph{Journal of the Operational Research Society}, 61:\penalty0
	523--529, 03 2010.
	
	\bibitem[Fleischmann(1985)]{BFleischmann_1985}
	B.~Fleischmann.
	\newblock A cutting plane procedure for the travelling salesman problem on road
	networks.
	\newblock \emph{European Journal of Operational Research}, 21\penalty0
	(3):\penalty0 307--317, 1985.
	
	\bibitem[Fukasawa et~al.(2018)Fukasawa, He, Santos, and
	Song]{RFukasawa_QHe_FSantos_YSong_2017}
	R.~Fukasawa, Q.~He, F.~Santos, and Y.~Song.
	\newblock A joint vehicle routing and speed optimization problem.
	\newblock \emph{INFORMS Journal on Computing}, 30:\penalty0 694--709, 2018.
	
	\bibitem[Garaix et~al.(2010)Garaix, Artigues, Feillet, and
	Josselin]{TGaraix_CArtigues_2010}
	T.~Garaix, C.~Artigues, D.~Feillet, and D.~Josselin.
	\newblock Vehicle routing problems with alternative paths: An application to
	on-demand transportation.
	\newblock \emph{European Journal of Operational Research}, 204:\penalty0
	62--75, 2010.
	
	\bibitem[García~Nájera and Bullinaria(2009)]{AGarciaNajera_JBullinaria_2009}
	A.~García~Nájera and J.~Bullinaria.
	\newblock Bi-objective optimization for the vehicle routing problem with time
	windows: Using route similarity to enhance performance.
	\newblock volume 5467, pages 275--289, 04 2009.
	
	\bibitem[Ghilas et~al.(2013)Ghilas, Demir, and {van Woensel}]{VGhilas_2013}
	V.~Ghilas, E.~Demir, and T.~{van Woensel}.
	\newblock Integrating passenger and freight transportation: model formulation
	and insights.
	\newblock BETA publicatie : working papers. Technische Universiteit Eindhoven,
	2013.
	
	\bibitem[Ghoseiri and Ghannadpour(2010)]{KGhoseiri_SGhannadpour_2010}
	K.~Ghoseiri and S.~Ghannadpour.
	\newblock Multi-objective vehicle routing problem with time windows using goal
	programming and genetic algorithm.
	\newblock \emph{Applied Soft Computing}, 10:\penalty0 1096--1107, 09 2010.
	
	\bibitem[Gupta et~al.(2015)Gupta, Ong, Zhang, Tan, Handa, Ishibuchi, and
	Tan]{AGupta_2015}
	A.~Gupta, Y.~Ong, A.~Zhang, P.~Tan, H.~Handa, H.~Ishibuchi, and K.C. Tan.
	\newblock A bi-level evolutionary algorithm for multi-objective vehicle routing
	problems with time window constraints.
	\newblock \emph{Proceedings of the 18th Asia Pacific symposium on intelligent
		and evolutionary systems—Volume 2. Proceedings in adaptation, learning and
		optimization}, 2:\penalty0 27--38, 01 2015.
	
	\bibitem[Huang et~al.(2006)Huang, Yao, and Raguraman]{BHuang_2006}
	B.~Huang, L.~Yao, and K.~Raguraman.
	\newblock Bi-level {GA} and {GIS} for multi-objective {TSP} route planning.
	\newblock \emph{Transportation Planning and Technology}, 29:\penalty0 105--124,
	2006.
	
	\bibitem[Karvonen et~al.(2018)Karvonen, Cugurullo, and
	Caprotti]{AKarvonen_FCugurullo_FCaprotti_2018}
	A.~Karvonen, F.~Cugurullo, and F.~Caprotti.
	\newblock \emph{Inside Smart Cities: Place, Politics and Urban Innovation}.
	\newblock Routledge, 2018.
	
	\bibitem[Kim et~al.(2015)Kim, Ong, Heng, Tan, and Zhang]{GKim_2015}
	G.~Kim, Y.~Ong, C.~Heng, P.~Tan, and A.~Zhang.
	\newblock City vehicle routing problem (city {VRP}): A review.
	\newblock \emph{IEEE Transactions on Intelligent Transportation Systems},
	16:\penalty0 1--13, 08 2015.
	
	\bibitem[Kramer et~al.(2014)Kramer, Subramanian, Vidal, and
	Cabral]{RKramer_2014}
	R.~Kramer, A.~Subramanian, T.~Vidal, and L.~Cabral.
	\newblock A matheuristic approach for the pollution-routing problem.
	\newblock \emph{European Journal of Operational Research}, 243, 04 2014.
	
	\bibitem[Kucharski and Drabicki(2017)]{RKucharski_ADrabicki_2021}
	R.~Kucharski and A.~Drabicki.
	\newblock Estimating macroscopic volume delay functions with the traffic
	density derived from measured speeds and flows.
	\newblock \emph{Journal of Advanced Transportation}, 2017, 2017.
	
	\bibitem[Kumar et~al.(2016)Kumar, Kondapaneni, Dixit, Goswami, Thakur, and
	Tiwari]{Kumar_Kondapaneni_Dixit_Goswami_Thakur_Tiwari_2016}
	R.S. Kumar, K.~Kondapaneni, V.~Dixit, A.~Goswami, L.S. Thakur, and M.K. Tiwari.
	\newblock Multi-objective modeling of production and pollution routing problem
	with time window: A self-learning particle swarm optimization approach.
	\newblock \emph{Computers \& Industrial Engineering}, 99:\penalty0 29 -- 40,
	2016.
	
	\bibitem[Letchford et~al.(2014)Letchford, Nasiri, and
	Oukil]{ALetchford_SNasiri_AOukil_2014}
	A.~Letchford, S.~Nasiri, and A.~Oukil.
	\newblock Pricing routines for vehicle routing with time windows on road
	networks.
	\newblock \emph{Computers \& Operations Research}, 51:\penalty0 331--337, 11
	2014.
	
	\bibitem[Li(2016)]{BLi_2016}
	B.~Li.
	\newblock Optimization of people and freight transportation : pickup and
	delivery problem variants.
	\newblock \emph{Ph.D. Thesis, TU/e Eindhoven University of Technology,
		Eindhoven, The Netherlands}, 2016.
	
	\bibitem[Ma and Xu(2014)]{YMa_2014}
	Y.~Ma and J.~Xu.
	\newblock Vehicle routing problem with multiple decision-makers for
	construction material transportation in a fuzzy random environment.
	\newblock \emph{International Journal of Civil Engineering}, 12:\penalty0
	332--346, 06 2014.
	
	\bibitem[Marinakis and Marinaki(2008)]{YMarinakis_2008}
	Y.~Marinakis and M.~Marinaki.
	\newblock A bilevel genetic algorithm for a real life location routing problem.
	\newblock \emph{International Journal of Logistics: Research and Applications},
	11:\penalty0 49--65, 02 2008.
	
	\bibitem[Marinakis et~al.(2007)Marinakis, Migdalas, and
	Pardalos]{YMarinakis_AMigdalas_PPardalos_2007}
	Y.~Marinakis, A.~Migdalas, and P.~Pardalos.
	\newblock A new bilevel formulation for the vehicle routing problem and a
	solution method using a genetic algorithm.
	\newblock \emph{Journal of Global Optimization}, 38:\penalty0 555--580, 08
	2007.
	
	\bibitem[Markevicius et~al.(2017)Markevicius, Navikas, Idzkowski, Valinevicius,
	Zilys, and Andriukaitis]{VMarkevicius_etal_2017}
	V.~Markevicius, D.~Navikas, A.~Idzkowski, A.~Valinevicius, M.~Zilys, and
	D.~Andriukaitis.
	\newblock Vehicle speed and length estimation using data from two anisotropic
	magneto-resistive (amr) sensors.
	\newblock \emph{Sensors}, 17:\penalty0 1--13, 08 2017.
	
	\bibitem[Miettinen(2012)]{KMiettinen_2012}
	K.~Miettinen.
	\newblock \emph{Nonlinear Multiobjective Optimization}, volume~12.
	\newblock Springer Science \& Business Media, New York, 2012.
	
	\bibitem[Mugayskikh et~al.(2018)Mugayskikh, Zakharov, and
	Tuovinen]{AVMugayskikh_VVZakharov_TTuovinen_2018}
	A.~V. Mugayskikh, V.~V. Zakharov, and T.~Tuovinen.
	\newblock Time- dependent multiple depot vehicle routing problem on megapolis
	network under wardrop's traffic flow assignment.
	\newblock \emph{2018 22nd Conference of Open Innovations Association (FRUCT)},
	pages 173--178, 2018.
	
	\bibitem[Nam and Pardo(2011)]{TNam_TAPardo_2011}
	T.~Nam and T.~A. Pardo.
	\newblock Conceptualizing smart city with dimensions of technology, people, and
	institutions.
	\newblock In \emph{Proceedings of the 12th Annual International Digital
		Government Research Conference: Digital Government Innovation in Challenging
		Times}, page 282–291, New York, NY, USA, 2011. Association for Computing
	Machinery.
	\newblock ISBN 9781450307628.
	
	\bibitem[Nasri et~al.(2018)Nasri, Bektaş, and Laporte]{MNasri_2018}
	M.~Nasri, T.~Bektaş, and G.~Laporte.
	\newblock Route and speed optimization for autonomous trucks.
	\newblock \emph{Computers \& Operations Research}, 100, 07 2018.
	
	\bibitem[Oyola et~al.(2016{\natexlab{a}})Oyola, Arntzen, and
	Woodruff]{JOyola_HArntzen_DWoodruff_2016a}
	J.~Oyola, H.~Arntzen, and D.~Woodruff.
	\newblock The stochastic vehicle routing problem, a literature review, part i:
	models.
	\newblock \emph{EURO Journal on Transportation and Logistics}, 7, 10
	2016{\natexlab{a}}.
	
	\bibitem[Oyola et~al.(2016{\natexlab{b}})Oyola, Arntzen, and
	Woodruff]{JOyola_HArntzen_DWoodruff_2016b}
	J.~Oyola, H.~Arntzen, and D.~Woodruff.
	\newblock The stochastic vehicle routing problem, a literature review, part ii:
	solution methods.
	\newblock \emph{EURO Journal on Transportation and Logistics}, 6, 10
	2016{\natexlab{b}}.
	
	\bibitem[Paszkowski et~al.(2021)Paszkowski, Herrmann, Richter, and
	Szarata]{JPaszkowski_etal_2021}
	J.~Paszkowski, M.~Herrmann, M.~Richter, and A.~Szarata.
	\newblock Modelling the effects of traffic-calming introduction to
	volume–delay functions and traffic assignment.
	\newblock \emph{Energies}, 14, 2021.
	
	\bibitem[Perron and Furnon()]{LPerron_VFurnon_2019}
	L.~Perron and V.~Furnon.
	\newblock Or-tools.
	\newblock URL \url{https://developers.google.com/optimization/}.
	
	\bibitem[Pillac et~al.(2013)Pillac, Gendreau, Guéret, and
	Medaglia]{VPillac_MGendreau_CGueret_AMedaglia_2013}
	V.~Pillac, M.~Gendreau, C.~Guéret, and A.L. Medaglia.
	\newblock A review of dynamic vehicle routing problems.
	\newblock \emph{European Journal of Operational Research}, 225:\penalty0 1 --
	11, 2013.
	
	\bibitem[Qian and Eglese(2014)]{JQian_REglese_2014}
	J.~Qian and R.~Eglese.
	\newblock Finding least fuel emission paths in a network with time-varying
	speeds.
	\newblock \emph{Networks}, 63:\penalty0 96--106, 2014.
	
	\bibitem[Raeesi and Zografos(2019)]{RRaeesi_KGZografos_2019}
	R.~Raeesi and K.~G. Zografos.
	\newblock The multi-objective steiner pollution-routing problem on congested
	urban road networks.
	\newblock \emph{Transportation Research Part B: Methodological}, 122:\penalty0
	457--485, 2019.
	
	\bibitem[Ritzinger et~al.(2016)Ritzinger, Puchinger, and
	Hartl]{URitzinger_JPuchinger_RHartl_2016}
	U.~Ritzinger, J.~Puchinger, and R.F. Hartl.
	\newblock A survey on dynamic and stochastic vehicle routing problems.
	\newblock \emph{International Journal of Production Research}, 54:\penalty0
	215--231, 2016.
	
	\bibitem[Rodríguez-Pereira et~al.(2019)Rodríguez-Pereira, Fernández,
	Laporte, Benavent, and Martínez-Sykora]{JRodriguezPereira_2019}
	J.~Rodríguez-Pereira, E.~Fernández, G.~Laporte, E.~Benavent, and
	A.~Martínez-Sykora.
	\newblock The steiner traveling salesman problem and its extensions.
	\newblock \emph{European Journal of Operational Research}, 278:\penalty0
	615--628, 2019.
	
	\bibitem[Sung and Nielsen(2020)]{ISung_PNielsen_2020}
	I.~Sung and P.~Nielsen.
	\newblock Speed optimization algorithm with routing to minimize fuel
	consumption under time-dependent travel conditions.
	\newblock \emph{Production \& Manufacturing Research}, 8:\penalty0 1--19, 2020.
	
	\bibitem[Tang et~al.(2019)Tang, Jayakar, Feng, Zhang, and
	Peng]{ZTang_KJayakar_XFeng_etal_2019}
	Z.~Tang, K.~Jayakar, X.~Feng, H.~Zhang, and R.~X. Peng.
	\newblock {Identifying smart city archetypes from the bottom up: A content
		analysis of municipal plans}.
	\newblock \emph{Telecommunications Policy}, 43\penalty0 (10), 2019.
	
	\bibitem[Vidal et~al.(2014)Vidal, Crainic, Gendreau, and
	Prins]{TVidal_TCrainic_MGendreau_CPrins_2014}
	T.~Vidal, T.G. Crainic, M.~Gendreau, and C.~Prins.
	\newblock A unified solution framework for multi-attribute vehicle routing
	problems.
	\newblock \emph{European Journal of Operational Research}, 234:\penalty0 658 --
	673, 2014.
	
	\bibitem[Vidal et~al.(2020)Vidal, Laporte, and Matl]{TVidal_2020}
	T.~Vidal, G.~Laporte, and P.~Matl.
	\newblock A concise guide to existing and emerging vehicle routing problem
	variants.
	\newblock \emph{European Journal of Operational Research}, 286:\penalty0 401 --
	416, 2020.
	
	\bibitem[W\"{a}chter and Biegler(2006)]{AWachter_LTBiegler_2006}
	A.~W\"{a}chter and L.~T. Biegler.
	\newblock On the implementation of an interior-point filter line-search
	algorithm for large-scale nonlinear programming.
	\newblock \emph{Math. Program.}, 106:\penalty0 25–57, mar 2006.
	
	\bibitem[Wenge et~al.(2014)Wenge, Zhang, Dave, Chao, and
	Hao]{RWenge_XZhang_CDave_LChao_SHao_2014}
	R.~Wenge, X.~Zhang, C.~Dave, L.~Chao, and S.~Hao.
	\newblock Smart city architecture: A technology guide for implementation and
	design challenges.
	\newblock \emph{China Communications}, 11:\penalty0 56--69, 2014.
	
	\bibitem[Xiong et~al.(2012)Xiong, Sheng, Rong, and Cooper]{ZXiong_2012}
	Z.~Xiong, H.~Sheng, W.~Rong, and D.~Cooper.
	\newblock Intelligent transportation systems for smart cities: A progress
	review.
	\newblock \emph{Science China Information Sciences}, 55, 12 2012.
	
	\bibitem[Yin et~al.(2015)Yin, Xiong, Chen, Wang, Cooper, and
	David]{CYin_ZXiong_HChen_2015}
	C.~Yin, Z.~Xiong, H.~Chen, J.~Wang, D.~Cooper, and B.~David.
	\newblock A literature survey on smart cities.
	\newblock \emph{Science China Information Sciences}, 58:\penalty0 1--18, 2015.
	
	\bibitem[Zografos and Androutsopoulos(2008)]{KGZografos_KNAndroutsopoulos_2008}
	K.~G. Zografos and K.~N. Androutsopoulos.
	\newblock A decision support system for integrated hazardous materials routing
	and emergency response decisions.
	\newblock \emph{Transportation Research Part C: Emerging Technologies},
	16:\penalty0 684--703, 2008.
	
\end{thebibliography}
\end{document}